\documentclass[a4paper,12pt]{amsart}
\usepackage[utf8]{inputenc}
\usepackage[english]{babel} 

\usepackage[T1]{fontenc}

\usepackage{amsthm}
\usepackage{amsmath}
\usepackage{amsfonts}
\usepackage{mathtools}
\usepackage[colorlinks]{hyperref}
\usepackage{enumerate}

\def\bdi{\begin{diagram}}
\def\edi{\end{diagram}}

\newtheorem{thm}{Theorem}[section]
\newtheorem{cor}[thm]{Corollary}
\newtheorem{lem}[thm]{Lemma}
\newtheorem{prop}[thm]{Proposition}
\theoremstyle{definition}
\newtheorem{defi}[thm]{Definition}
\newtheorem{defis}[thm]{Definitions}
\newtheorem{conj}[thm]{Conjecture}
\newtheorem{conv}[thm]{Convention}
\newtheorem{nota}[thm]{Notation}
\newtheorem{rem}[thm]{Remark}
\newtheorem{rems}[thm]{Remarks}
\newtheorem{exa}[thm]{Example}
\newtheorem{exas}[thm]{Examples}
\newtheorem{prob}[thm]{Problem}
\newtheorem{probs}[thm]{Problems}
\newtheorem{ques}[thm]{Question}
\newtheorem{sett}[thm]{Setting}
\newtheorem{sit}[thm]{}

\newcommand{\sing}{\operatorname{{\rm sing}}}

\newcommand{\Ker}{ \operatorname{{\rm Ker}}}

\newcommand{\ED}{ \operatorname{{\rm ED}}}

\newcommand{\Aut}{ \operatorname{{\rm Aut}}}

\newcommand{\LND}{\operatorname{{\rm LND}}}
\newcommand{\GL}{\operatorname{{\rm GL}}}
\newcommand{\SL}{ \operatorname{{\rm SL}}}

\newcommand{\Hom}{ \operatorname{{\rm Hom}}}

\def\pr{\mathop{\rm pr}}

\def\reg{{\mathop{\rm reg}}}
\def\sing{{\mathop{\rm sing}}}

\def\codim{\mathop{\rm codim}}

\def\Pic{\mathop{\rm Pic}}

\def\lto{\longrightarrow}

\renewcommand{\epsilon}{\varepsilon}

\def\and{\quad\mbox{and}\quad}

\newcommand{\C}{\ensuremath{\mathbb{C}}}
\newcommand{\T}{\ensuremath{\mathbb{T}}}
\newcommand{\R}{\ensuremath{\mathbb{R}}}

\newcommand{\Z}{\ensuremath{\mathbb{Z}}}

\newcommand{\G}{\ensuremath{\mathbb{G}}}

\newcommand{\bk}{{\ensuremath{\rm \bf k}}}

\newcommand{\tC}{{\tilde C}}
\newcommand{\tD}{{\tilde D}}
\newcommand{\tF}{{\tilde F}}

\newcommand{\tH}{{\tilde H}}

\newcommand{\tV}{{\tilde V}}

\newcommand{\tZ}{{\tilde Z}}

\newcommand{\bZ}{{\bar Z}}

\newcommand{\tO}{{\tilde O}}
\newcommand{\tT}{{\tilde T}}

\newcommand{\cW}{{\ensuremath{\mathcal{W}}}}

\newcommand{\cG}{{\ensuremath{\mathcal{G}}}}
\newcommand{\cS}{{\ensuremath{\mathcal{S}}}}

\newcommand{\cA}{{\ensuremath{\mathcal{A}}}}

\newcommand{\cO}{{\ensuremath{\mathcal{O}}}}
\newcommand{\cC}{{\ensuremath{\mathcal{C}}}}

\newcommand{\cN}{{\ensuremath{\mathcal{N}}}}

\newcommand{\cX}{{\ensuremath{\mathcal{X}}}}
\newcommand{\cY}{{\ensuremath{\mathcal{Y}}}}

\newcommand{\cZ}{{\ensuremath{\mathcal{Z}}}}
\newcommand{\p}{\partial}

\renewcommand{\rho}{\varrho}

\def\bals#1\eals{\begin{align*}#1\end{align*}}
\def\bal#1\eal{\begin{align}#1\end{align}}

\def\SAut{\mathop{\rm SAut}}
\def\SL{\mathop{\rm SL}}

\def\A{{\mathbb A}}

\def\PP{{\mathbb P}}

\def\dd{ {\rm d }}

\renewcommand{\phi}{\varphi}

\newcommand{\bnum}{\begin{enumerate}}
\newcommand{\enum}{\end{enumerate}}

\newcommand{\brem}{\begin{rem}}
\newcommand{\brems}{\begin{rems}}
\newcommand{\erem}{\end{rem}}
\newcommand{\erems}{\end{rems}}
\newcommand{\bprob}{\begin{prob}}
\newcommand{\eprob}{\end{prob}}
\newcommand{\bprobs}{\begin{probs}}
\newcommand{\eprobs}{\end{probs}}
\newcommand{\bques}{\begin{ques}}
\newcommand{\eques}{\end{ques}}
\newcommand{\bexa}{\begin{exa}}
\newcommand{\bexas}{\begin{exas}}
\newcommand{\eexa}{\end{exa}}
\newcommand{\eexas}{\end{exas}}
\newcommand{\bdefi}{\begin{defi}}
\newcommand{\edefi}{\end{defi}}
\newcommand{\bdefis}{\begin{defis}}
\newcommand{\edefis}{\end{defis}}
\newcommand{\bcor}{\begin{cor}}
\newcommand{\ecor}{\end{cor}}
\newcommand{\blem}{\begin{lem}}
\newcommand{\elem}{\end{lem}}
\newcommand{\bconv}{\begin{conv}}
\newcommand{\econv}{\end{conv}}
\newcommand{\bconj}{\begin{conj}}
\newcommand{\econj}{\end{conj}}
\newcommand{\bprop}{\begin{prop}}
\newcommand{\eprop}{\end{prop}}
\newcommand{\bthm}{\begin{thm}}
\newcommand{\ethm}{\end{thm}}
\newcommand{\bnota}{\begin{nota}}
\newcommand{\enota}{\end{nota}}
\newcommand{\bsit}{\begin{sit}}
\newcommand{\esit}{\end{sit}}
\newcommand{\be}{\begin{equation}}
\newcommand{\ee}{\end{equation}}
\newcommand{\bproof}{\begin{proof}}
\newcommand{\eproof}{\end{proof}}
\newcommand{\bsett}{\begin{sett}}
\newcommand{\esett}{\end{sett}}
\def\ba{\begin{array}}
\def\ea{\end{array}}

\newcommand{\sslash}{\mkern-3.5mu\mathbin{/\mkern-6mu/}\mkern-3.5mu}

\begin{document}

\title{Lines in affine toric varieties}

\author{Shulim Kaliman}
\address{
University of Miami, Department of Mathematics, Coral Gables, FL 33124, USA}
\email{kaliman@math.miami.edu}

\date{\today}
\maketitle

\begin{abstract}
We prove that up to automorphisms of the target the affine line $\A^1$ admits a unique embedding into the regular part of
an affine simplicial toric variety of dimension at least 4 which is smooth in codimension 2.
This is an analog of the well-known result on the existence of a linearization of any 
polynomial embedding $\A^1 \hookrightarrow \A^n$ for $n\geq 4$.
\end{abstract}



\section{Introduction}

Let $\varphi : C\to C'$ be an isomorphism of two smooth polynomial curves contained in the regular part $Y_{\rm reg}$ of an affine algebraic variety $Y$
over an algebraically closed field $\bk$ of characteristic zero.
It may happen that $\varphi$ extends to an automorphism of $Y$ and our first aim is to describe some affine algebraic varieties for
which this extension takes place. 

This problem was studied in several papers \cite{AbMo}, \cite{Su}, \cite{Cr},\cite{Je}, \cite{St},\cite{FS}, \cite{Ka20}
and \cite{AZ}. 
It turns out that the answer is positive for some classes of flexible varieties of dimension $n\geq 4$ where $Y$ is flexible if the subgroup $\SAut (Y)$ of the automorphism
group $\Aut (Y)$ of $Y$ generated by all one-parameter unipotent subgroups acts transitively on $Y_{\rm reg}$. Say, this is so if $Y=\A^n$ with $n\geq 4$ \cite{Cr},\cite{Je}.
For $n=3$ the answer is unknown but for $n=2$ the famous Abhyankar-Moh-Suzuki theorem \cite{AbMo},\cite{Su} states that an isomorphism of two smooth plane polynomial
curves always extends to an automorphism of the plane $\A^2$. Perhaps, $\A^2$ is the only example of a two-dimensional flexible variety
with this property. If $Y$ is an affine simplicial toric variety $\A^2/G$ where $G$ is a finite subgroup of $\SL_2(\bk)$ acting naturally on $\A^2$,
then Arzhantsev and Zaidenberg \cite{AZ} showed that the answer is negative. They actually classified  up to automorphisms of $Y$ all  smooth polynomial curves in $Y_{\rm reg}$
(there are only finite number of isomorphism classes of such curves). 

In this paper we study the case when $Y$ is an affine simplicial toric variety of dimension $n \geq 4$
(i.e., $Y=\A^n/G$ where $G$ is a finite subgroup of $\SL_n(\bk)$ acting naturally on $\A^n$). We show that the answer to this extension problem is
positive under the assumption of smoothness in codimension 2. 

Furthermore, recall that given a subvariety $Z$ of $Y$ with defining ideal $I$ in the algebra $\bk [Y]$ of regular functions on $Y$
its $k$th infinitesimal neighborhood is the scheme with the defining ideal $I^k$. In particular, if $W$ is another subvariety of $Y$ with
defining ideal $J$, then an isomorphism $\cZ \to \cW$  of $k$th infinitesimal neighborhoods of $Z$ and $W$ is determined by an isomorphism of algebras
$\frac{\bk [Y]}{I^k}\to \frac{\bk [Y]}{J^k}$.  There are natural obstacles for extending such isomorphisms to automorphisms of $Y$.
Say,  let $Z=W$ be a strict complete intersection given in $Y$ by $u_1=\ldots =u_{m}=0$.
Then an automorphism $\psi : \frac{\bk [Y]}{I^k}\to \frac{\bk [Y]}{I^k}$ over $\frac{\bk [Y]}{I}$ is given by polynomials 
$f_1, \ldots ,f_{ m}$ in $u_1, \ldots, u_{ m}$
over $\frac{\bk [Y]}{I}$ of degree at most $k-1$. If $Y$ does not admit nonconstant invertible functions and
$\psi$ is extendable to an automorphism of $Y$, then one can see that the Jacobian 
$\det\left[ \frac{\p f_i}{\p u_j}\right]_{i,j=1}^{ m}$ must be equal to a
nonzero constant  modulo $I^{k-1}$ in which case we say that $\psi$ has a nonzero constant Jacobian. There is also a notion of a nonzero constant Jacobian of 
an isomorphism $\cZ \to \cW$ in the case when both $Z$ and $W$ are smooth polynomial curves in a normal toric variety $Y$ contained in $Y_{\rm reg}$ (see Definition \ref{sim.d1}).
 The question when such isomorphisms with nonzero constant Jacobians are extendable to automorphisms of $Y$
was considered in \cite{KaUd} and \cite{Ud}. In combinations with the results of \cite{KaUd} and \cite{Ud} we get our first main result (Corollary \ref{sim.c1}).

\bthm\label{int.t1} Let $Y$ be  an affine simplicial toric variety  smooth in codimension 2 such that
 $\dim Y \geq 4$. Let $\varphi : \cC_1 \to \cC_2$ be an isomorphism
of $k$th infinitesimal neighborhoods of two smooth polynomial curves contained in  $Y_{\rm reg}$ such that the Jacobian of
$\varphi$ is a nonzero constant.        
Then $\varphi$ extends to an automorphism of $Y$.
\ethm

The second subject of this paper is related to the theorem of Holme \cite[Theorem 7.4]{Hol} (later rediscovered in \cite{Ka91} and \cite{Sr}).
It states that if $Z$ is an affine algebraic variety with $$\ED(Z):=\max (2\dim Z+1, \dim TZ) \leq n,$$ then $Z$ admits a closed embedding into $\A^n$
(the version of this theorem with a smooth $Z$ appeared originally in  \cite[Theorem 2.1]{Swan}).
Recently Feller and van Santen \cite{FS21}
proved that if $X$ is an affine algebraic variety isomorphic to a simple linear algebraic group and $Z$ is smooth, then $Z$ admits a closed
embedding into $X$, provided that $\dim X>\ED (Z)=2\dim Z+1$.  Since affine spaces, simple linear algebraic groups and normal affine toric varieties are examples of flexible varieties
it is natural to look for analogues of Holme's theorem in the flexible case. In this paper we prove the  following.

\bthm\label{int.t2} {\rm (Theorem \ref{fle.t2})}  Let $Z$ be an affine algebraic variety and $X$ be a smooth quasi-affine flexible variety of dimension at least $\ED (Z)$.
Then $Z$ admits an injective immersion into $X$.
\ethm

In the case when $X$ is a normal affine toric variety we also find conditions which guarantee that $Z$ admits a closed embedding into $X$ (Theorem \ref{emb.t1}).
The formulation of the latter theorem is subtler when $X$ is simplicial and it is a consequence of the following more general fact.

\bthm\label{int.t3} {\rm (Corollary \ref{fle.c1})}
Let $\psi : \A^n \to Y$ be a finite morphism where $Y$ is normal. Suppose that $Z$ is an affine algebraic variety
such that $\ED ( Z) \leq n$ and  $\dim Z < \codim_YY_\sing$. Then $Z$ admits 
a closed embedding in $Y$ with the image contained in $Y_{\rm reg}$.   
\ethm

The paper is organized as follows. In Section 2 we survey  the technique developed in \cite{Ka20} which was later clarified in \cite{KaUd}.
In particular, one can find there formal definitions of
locally nilpotent vector fields and flexible varieties. Section 2 contains also a modified version of Theorem 4.2 from \cite{Ka20}
which is a crucial tool in this paper.
Using this result we prove Theorems \ref{int.t2} and \ref{int.t3} in Section 3.
In Section 4 we introduce notations for toric varieties which are used freely throughout the rest of the paper
and prove some simple facts about normal affine toric varieties. Section 4 contains an analogue of Holme's theorem for normal affine toric varieties.
In Section 5 we study locally nilpotent vector fields on normal affine toric varieties with no torus factors. The properties of locally nilpotent vector fields are crucial for
us since compositions of elements of the flows of such vector fields produce automorphisms that extend isomorphisms of smooth polynomial curves.
In Section 6 we prove Theorem \ref{int.t1}.\\

{\em Acknowlegement.}  The author is grateful to the referee for very useful comments and corrections.

\section{Flexible varieties: preliminaries}

In this section we present some technical tools developed in \cite{Ka20} with later clarifications in \cite{KaUd} which we use in this paper.
We shall also give a modified version of  \cite[Theorem 4.2]{Ka20}.

\bdefi\label{def.d1}  (1) Given an irreducible algebraic variety $\cA$ and
a map $\varphi:\cA\to\Aut(X)$  (where $\Aut (X)$ is the group of algebraic automorphisms of $X$) we say that $(\cA,\phi)$
is an {\em algebraic family of automorphisms of $X$} if the induced map
$\cA\times X\to X$, $(\alpha,x)\mapsto \varphi(\alpha).x$ is a morphism (see \cite{Ra}).

(2)  If we want to emphasize additionally that  $\varphi (\cA)$ is contained in a subgroup $G$ of $\Aut (X)$, then we say that
$\cA$ is an {\em algebraic  $G$-family} of automorphisms of $X$.            

(3) In the case when $\cA$ is a connected algebraic group and the induced map 
$\cA\times X\to X$ is not only a morphism but also an action of $\cA$ on $X$ we call this family a {\em connected algebraic subgroup} of $\Aut (X)$.

(4)
Following \cite[Definition 1.1]{AFKKZ} we call a subgroup $G$ of $\Aut (X)$ algebraically generated if it is generated as an abstract group by a family 
$\cG$ of connected algebraic subgroups of $\Aut (X)$.
\edefi

We have the following important fact \cite[Theorem 1.15]{AFKKZ} (which is the analogue of  the Kleiman transversality theorem \cite{Kl}
for algebraically generated groups).

\bthm\label{agga.t1} ({\rm Transversality Theorem}) Let a subgroup $G\subseteq \Aut(X)$ be
algebraically generated by a system $\cG$ of connected algebraic
subgroups closed under conjugation in $G$. Suppose that $G$ acts
with an open orbit $O\subseteq X$.

Then there exist subgroups $H_1,\ldots, H_m\in \cG$ such that for
any locally closed reduced subschemes $Y$ and $Z$ in $O$ one can
find a Zariski dense open subset $U=U(Y,Z)\subseteq H_1\times
\ldots \times H_m$ such that every element $(h_1,\ldots, h_m)\in
U$ satisfies the following: \bnum[(a)]
\item  {\em The translate $(h_1\cdot\ldots\cdot h_m).Z_\reg$
meets $Y_\reg$
transversally. } 
\item $\dim (Y\cap (h_1\cdot\ldots\cdot h_m).Z)\le
\dim Y+\dim Z-\dim X$. \footnote{We put the dimension of empty sets equal to $-\infty$.}\\
In particular $Y\cap (h_1\cdot\ldots\cdot h_m).Z=\emptyset$
if $\dim Y+\dim Z<\dim X$.
\enum
\ethm

\bdefi\label{def.d3}  (1) A nonzero derivation $\delta$ on the ring $A$ of regular functions on an affine algebraic variety $X$ is called {\em  locally nilpotent}
if for every $0\ne a \in A$ there exists a natural $n$ for which $\delta^n (a)=0$. 
This derivation can be viewed as a vector field on $X$ which
we also call {\em locally nilpotent}. The set of all locally nilpotent vector fields on $X$ will be denoted by $\LND (X)$. 
The flow of $\delta \in \LND (X)$ 
 is an algebraic $\G_a$-action on $X$, i.e., the action of the group $(\bk, +)$ 
which can be viewed as a one-parameter unipotent group $U$ in the group $\Aut (X)$ of all algebraic automorphisms of $X$.
In fact, every $\G_a$-action is a flow of a locally nilpotent vector field (e.g, see \cite[Proposition 1.28]{Fre}).

(2) If $X$ is a quasi-affine variety, then an algebraic vector field $\delta$ on $X$ is called locally nilpotent if $\delta$ extends
to a locally nilpotent vector field  $\tilde \delta$ on some affine algebraic variety $Y$ containing $X$ such
that $\tilde \delta$ vanishes on $Y\setminus X$ where ${\rm codim}_Y (Y \setminus X) \geq 2$.  Note that under this assumption
$\delta$ generates a $\G_a$-action on $X$ and we use again the notation $\LND (X)$ for the set of all locally nilpotent vector
fields on $X$.\edefi

\bdefi\label{def.d4}

(1) For every locally nilpotent vector fields $\delta$ and each function $f \in \Ker \delta$ from its kernel the field
$f\delta$ is called a replica of $\delta$. Recall that such replica is automatically locally nilpotent.

(2)  Let $\cN$ be a set of locally nilpotent vector fields on $X$ and $G_\cN \subset  \Aut (X)$ denote the group generated
by all flows of elements of $\cN$. We say that $G_\cN$ {\em is generated by $\cN$}.

(3) A collection of locally nilpotent vector fields $\cN$ is called saturated if $\cN$ is closed under conjugation by elements in $G_\cN$
and for every $\delta \in \cN$ each replica of $\delta$ is
also contained in $\cN$.

\edefi

\bdefi\label{def.d5}  Let $X$ be a normal quasi-affine algebraic variety of dimension at least 2,
$\cN$ be a saturated set of locally nilpotent vector fields on $X$ and  $G=G_\cN$ be the group
generated by  $\cN$.
Then $X$ is called $G$-flexible if for any point $x$ in the smooth part $X_{\rm reg}$ of $X$ the vector space $T_xX$ is generated 
by  the values of locally nilpotent vector fields from $\cN$ at $x$
(which is equivalent to the fact that $G$ acts transitively on $X_{\rm reg}$ \cite[Theorem 2.12]{FKZ}). In the case of $G=\SAut (X)$ we call $X$ flexible
without referring to $\SAut (X)$ (recall that $\SAut (X)$ is the subgroup of $\Aut X$ generated by all one-parameter unipotent subgroups).
\edefi

The following is a modified version of \cite[Theorem 4.2]{Ka20}.

\bthm\label{gp1.t1}  
Let $X$ be a smooth algebraic variety, $Q$ be a normal algebraic variety, $\rho : X \to Q$ be a dominant morphism
and $G \subset \Aut (X)$ be an algebraically generated group acting 2-transitively on  $X$.
 Suppose that  $Q_0$ is a smooth open dense subset of $Q$,  $X_0=\rho^{-1} (Q_0)$
and $Z$ is a locally closed reduced subvariety of $X$.

{\rm (i)}  Suppose that
\be\label{fle.eq1} \dim X_0\times_{Q_0} X_0= 2\dim X -\dim Q.\ee
and $\dim Q \geq \dim Z+m$ where $m \geq 1$. Then  there exists an algebraic $G$-family $\cA$ of automorphisms
of $X$  such that
for a general element $\alpha \in \cA$ one can find 
a constructible subset $R$ of $\alpha (Z)\cap X_0$ of dimension $\dim R \leq \dim Z-m$ for which
$\rho (R)$ and $ \rho (\alpha (Z)\setminus R)$ are disjoint  
and the restriction 
$$\rho|_{(\alpha (Z)\cap X_0)\setminus R}:  (\alpha (Z) \cap X_0) \setminus R\to  Q_0$$
of $\rho$ is injective.
In particular, if $\dim Q\geq  2\dim Z +1$  and  $Z_\alpha'=\rho \circ \alpha (Z)$, then
for a general element $\alpha \in \cA$ the morphism
$$\rho|_{\alpha (Z)\cap X_0}:  \alpha (Z) \cap X_0 \to Z_\alpha'\cap Q_0$$ 
is a bijection,
while in the case of a pure-dimensional $Z$ and $\dim Q \geq \dim Z +1$ 
 this morphism   is birational.

{\rm (ii)} Let $G$ be generated by a saturated set $\cN$ of locally nilpotent vector fields on $X$
(in particular, $X$ is $G$-flexible) and  $Y=\bigcup_{x \in X_0} \Ker \{ \rho_* : T_x X_0 \to T_{\rho (x)} Q_0\}$. 
Let
\be\label{fle.eq2}      \dim Y =\dim TX -\dim Q.            \ee
Suppose that $\dim TZ  \leq \dim Q$. 
Then there exists an algebraic family $\cA$ of $G$-automorphisms of $X$ such that
for a general element $\alpha \in \cA$  and
every $z \in \alpha (Z) \cap X_0$ the induced map $\rho_* : T_z\alpha (Z) \to T_{\rho (z)}Q$ of the tangent spaces is injective.

\ethm

\bproof  For every variety $\cX$  denote by $S_\cX $ the variety $S_\cX=(\cX\times \cX) \setminus \Delta_\cX$ where $\Delta_\cX$ is the diagonal in $\cX \times \cX$. 
Then every automorphism in $\Aut (X)$ can be lifted to an automorphism of $S_X$. In particular, we have a $G$-action on $S_X$ and by the assumption this action is transitive
on $S_X$.
Consider the subvariety  $Y=(X_0\times_{Q_0} X_0) \setminus \Delta_X \subset S_X$.
By Formula \eqref{fle.eq1} $\dim Y=2\dim X -\dim Q$ 
(i.e., the codimension of $Y$ in $S_X$ is $\dim Q$). 
By Theorem \ref{agga.t1} (b) we can choose algebraic subgroups $H_1, \ldots, H_m$ of $G$ such that for a general element $(h_1, \ldots , h_m) \in H_1 \times \cdots \times H_m$
one has $$\dim W \leq \dim Y+  \dim S_Z -\dim S_X= \dim S_Z-\dim Q \leq 2 \dim Z-\dim Q$$ where 
$W= Y \cap \alpha (S_Z)$ for $\alpha =h_1 \cdot \ldots \cdot h_m$.
Hence, in case (i) the dimension of $W$ is at most $\dim Z-m$. Let $R$ be the image of $W$ under one of the two natural projections $X\times_QX \to X$. 
In particular, $R$ is a constructible set by  Chevalley's theorem \cite[Chap. II, Exercise 3.19]{Har}, $R\subset \alpha (Z)\cap X_0$ and  $\dim R \leq \dim Z-m$.
Note that for $z \in \alpha (Z)\cap X_0$ one has $\rho^{-1}(\rho (z))\cap \alpha (Z)=z$ iff $z \notin R$.
Hence, the restriction of $\rho$ to $(\alpha (Z)\cap X_0)\setminus R$ is injective.
Therefore, letting $\cA = H_1 \times \cdots \times H_m$, we get (i).

In (ii)  for every variety $\cX$ and a subvariety $\cY$ of the tangent bundle $T\cX$ let $\cY^*=\cY \setminus \cS $ where $\cS$
is the zero section of the natural morphism $T\cX \to \cX$.
Every automorphism $\alpha \in \Aut (X)$ generates an automorphism of $TX$. In particular, $G$ acts on $(TX)^*$
and by  \cite[Theorem 4.11 and Remark 4.16]{AFKKZ} this action is transitive.
By Formula \eqref{fle.eq2} $\dim Y^*=\dim TX -\dim Q$. 
By Theorem \ref{agga.t1} 
 we can choose one-parameter unipotent algebraic subgroups $\tH_1, \ldots, \tH_{\tilde m}$ of $G$ such that 
for a general element $(\tilde h_1, \ldots , \tilde h_{\tilde m}) \in \tH_1 \times \cdots \times \tH_{\tilde m}$ and 
$Z''=(\tilde h_1 \cdot \ldots \cdot \tilde h_{\tilde m}) (Z)$ one has $\dim Y^* \cap (TZ'')^* \leq \dim Y^* + \dim (TZ)^* -\dim TX^*\leq 0$.
Note that if $Y^* \cap (TZ'')^*$ contains a point, then $\dim Y^* \cap TZ''$ must be at least 1
 (since this point is a vector in $TZ''$ and then $Y^* \cap (T(Z'')^*$ 
contains all nonzero vectors proportional to that one).
That is, $Y^*\cap (T(Z'')^*=\emptyset$.
This implies that for every $z\in Z''\cap X_0$ the restriction of $\rho_*$ to 
$T_zZ''$ is injective. Consequently, the restriction of $\rho_*$ to 
$T_zZ''$ is injective i.e., we have (ii).
 \eproof

Let us describe some $G$-families $\cA$ satisfying the conclusions of Theorem \ref{gp1.t1}.

\bdefi\label{prep.d6} Let $X$ be a smooth algebraic variety and $G$ be a subgroup of $X$.  
Consider $(X\times X)\setminus \Delta$ (where $\Delta$ is the diagonal), the complement $(TX)^*$ to the zero section in the tangent bundle
of $X$ and the frame bundle ${\rm Fr} (X)$ of $TX$ (i.e., the fiber of  ${\rm Fr} (X)$ over $x \in X$ consists of all bases of $T_xX$).
Projectivization of $TX$ replaces ${\rm Fr}(X)$ with a bundle  ${\rm PFr} (X)$ whose fiber over $x$ consists of all ordered $n$-tuples of points
in the proectivization $\PP^n$ of $T_xX$  (where $n=\dim X$) that are not contained in the same hyperplane of $\PP^n$. 
Then we have natural $G$-actions on all these objects.
Let $Y$ be either $X$, or $(X\times_PX)\setminus \Delta$, or $(TX)^*$,  or ${\rm PFr} (X)$. Suppose that  the $G$-action is transitive on $Y$.
Then we say that an algebraic $G$-family $\cA$ of automorphisms of $X$ is
{\em a regular $G$-family  for $Y$} if

(i) $\cA= H_{m}\times\ldots \times H_{1}$ where each $H_i$ belongs to $\cG$;

(ii) for a suitable open dense subset $U\subseteq H_{m}\times\ldots \times H_{1}$, the map
\be\label{prep.eq2}\begin{array}{c}  \Psi: H_m\times \ldots\times H_1\times Y \lto Y\times Y,\\
 (h_m,\ldots,h_1,y)\mapsto
((h_m\cdot\ldots\cdot h_1).y ,y) \end{array}
\ee
is smooth on $U\times Y$. 

An algebraic $G$-families $\cA$ that are regular for all four varieties
$X$,  $(X\times_PX)\setminus \Delta$,  $(TX)^*$ and ${\rm PFr} (X)$ will be  called  a {\em perfect $G$-family for $Y$}.
\edefi

\bprop\label{new.p1} Let $X$ be a smooth algebraic variety and $G\subset \Aut (X)$ be a group algebraically generated by
a family $\cG$ of algebraic connected subgroups of $\Aut (X)$. Suppose that $G$ acts transitively on $X$.

{\rm (1)} Then there exists a regular $G$-family for $X$ (which is of the form $\cA =H_1 \times \ldots \times H_m$ where each $H_i$
is an element of $\cG$).

{\rm (2)} Every regular $G$-family for $X$ satisfies the conclusions of Theorem \ref{agga.t1}.

{\rm (3)} If $\cA$ is a regular (resp. perfect) $G$-family  for $X$ and $H$ is an element of $\cG$ then $H\times \cA$ and $\cA \times H$
are also regular (resp. perfect) $G$-families for $X$.

{\rm (4)} In particular, if $X$ is $G$-flexible, then there exists a perfect $G$-family.

{\rm (5)} Let $X$ be $G$-flexible. Every $G$-family regular for  $X$ (resp. for $(TX)^*$) 
satisfies the conclusion of Theorem \ref{gp1.t1} (i) (resp. (ii)). In particular,
every perfect $G$-family satisfies the conclusion  of Theorem \ref{gp1.t1} (i)-(ii).

\eprop

\bproof Statement (1) is proven in \cite[Proposition 1.15]{AFKKZ}. Statements (2) and (3) are proven in \cite[Proposition 1.10]{Ka20}.
The fourth statement follows from the fact that in the flexible case for every $m>0$ the group $G$ acts $m$-transitively on $X$ and
also transitively on $(TX)^*$ and ${\rm PFr} (X)$ \cite[Theorem 4.11 and
Remark 4.16]{AFKKZ}. Modulo the definitions of regular and perfect families  and (2) statement (5) follows from the proof  of Theorem \ref{gp1.t1}.
\eproof

\bprop\label{fle.p1}
Let $X$ be a smooth algebraic variety, $Q$ be a normal algebraic variety and $\rho : X \to Q$ be a dominant morphism.
Let $Q_0$ be  a smooth dense open subset   of $Q$ and $X_0=\rho^{-1} (Q_0)$.
Suppose that for every $q \in Q_0$ the fiber $\rho^{-1} (q)$ is smooth and of dimension $\dim X -\dim Q$.
Then Formulas \eqref{fle.eq1} and \eqref{fle.eq2}
hold. In particular, if $X$ is $G$-flexible and $Z$ satisfies the assumption of Theorem \ref{gp1.t1}(i)-(ii),
then for every perfect $G$-family  $\cA$ of automorphisms of $X$ and a general $\alpha\in \cA$
the morphism $\rho|_{\alpha (Z) \cap X_0} : \alpha (Z) \cap X_0 \to Q_0$ is an injective immersion. 
\eprop

\bproof The validity of Formula \eqref{fle.eq1} is straightforward. 
 Consider any irreducible subvariety
 $P$ of $Q_0$ and an  irreducible component $W$ of $\rho^{-1}(P)$ whose image is dense in $P$.
 Since every fiber of $\rho$ is smooth so is the generic fiber of $\rho|_W$ (e.g., see \cite[Lemma 2.1]{KrRu}).
 Hence, replacing  $P$ by its open dense subset we can suppose that $W$ is smooth.                      
Furthermore,  by \cite[Chapter III, Corollary 10.7]{Har} 
 we can suppose that  $\rho|_W: W \to P$ is smooth. 
Therefore, $\dim  \Ker \rho_* |_{ T_x X}\leq \dim X -\dim P$ for every $x \in W$. 
This implies that the dimension of $ Y_W=\bigcup_{x \in W} \Ker \{ \rho_* : T_x X \to T_{\rho (x)} Q\}$
is at most $\dim X -\dim P+ \dim W$. Since $\rho|_{X_0}$ is equidimensional one has $\dim W -\dim P =\dim X -\dim Q$
and $\dim Y_W \leq \dim TX-\dim Q$. Of course, we can suppose that the latter inequality is true for every irreducible
component $W$ in $\rho^{-1} (P)$ which yields the desired conclusion. 
\eproof

\section{Embedding theorems for flexible varieties}

\bnota\label{emb.d1} Let $Z$ be an affine algebraic variety and $TZ$ be its Zariski tangent bundle. Then
we let $$\ED (Z)=\max (2\dim Z+1, \dim TZ).$$
\enota

By \cite[Theorem 7.4]{Hol} for every affine algebraic variety $Z$
there exists a closed embedding of $Z$ into $\A^{\ED (Z)}$.

\bthm\label{emb.t2} Let $\psi : W \to Y$ be a finite morphism where $W$ is a smooth  flexible variety and $Y$ is normal. Let $Z$ be a quasi-affine algebraic variety
which admits a closed embedding in $W$. Suppose also that $\dim Z < \codim_Y Y_{\sing}$. Then $Z$ admits 
a closed embedding in $Y$ with the image contained in $Y_{\rm reg}$. 
\ethm

\bproof  One can treat $Z$ as a closed subvariety of $W$. By  Theorem \ref{agga.t1}  there exists an algebraic family $\cA$
of automorphisms of $W$ such that for a general $\alpha\in \cA$ the variety $\alpha (Z)$ does not meet $\psi^{-1} (Y_\sing)$.
By Proposition \ref{new.p1}(2)-(4) enlarging $\cA$ we can suppose that it is a perfect family.
Proposition \ref{fle.p1}  implies now
that $\psi|_{\alpha (Z)} : \alpha (Z) \to Y_\reg \subset Y$ is an injective immersion. Since $\psi$ is finite $\psi|_{\alpha (Z)}$ is also proper.  Hence, we are done.
\eproof

\bcor\label{fle.c1} 
Let $\psi : \A^r \to Y$ be a finite morphism where $Y$ is normal. Suppose that $Z$ is an affine algebraic variety
such that $\ED ( Z) \leq r$ and  $\dim Z < \codim_YY_\sing$. Then $Z$ admits 
a closed embedding in $Y$ with the image contained in $Y_{\rm reg}$.   
\ecor

\brem\label{fle.r1} For every $m>0$ there are examples of affine algebraic varieties of dimension $m$ that cannot be embedded in $\A^{2m}$ \cite{BMS}. In particular, 
Holme's theorem is sharp and we cannot improve
the assumption $\ED ( Z) \leq r$ in Corollary \ref{fle.c1}. However, the author does not know if 
the assumption $\dim Z < \codim_YY_\sing$ is optimal for every $Y$ as in Corollary \ref{fle.c1} (especially, in the light of Theorem \ref{fle.t2} below).
\erem

\bprop\label{pre.p2} Let $X$ be a $G$-flexible variety and $H= H_m \times \ldots \times H_1$ be a perfect $G$-family of automorphisms of
$X$ (where $H_1, \ldots, H_m$ are unipotent subgroups of $G$). 
Suppose that an open  dense supset $U \subset  H$ is such that the morphism $\Psi : H \times Y \to Y \times Y$ as in Formula \eqref{prep.eq2}
is smooth on $U\times Y$ for every $Y$ equal to one of the varieties $X$,  $(X\times_PX)\setminus \Delta$,  $(TX)^*$ and ${\rm PFr} (X)$.
Then $m, H_1 \ldots, H_m$ and $U$ can be chosen so that the codimension of $H\setminus U$ in $H$ is arbitrarily large.
\eprop

\bproof  Since one can increase $m$ by Proposition \ref{new.p1}(3)  the desired conclusion
follows from \cite[page 778, footnote]{AFKKZ}. \eproof

\bprop\label{pre.p3} Let $H$ be a smooth flexible variety, $X$ be a normal algebraic variety and $\varphi : H \to X$ be a dominant morphism
such that $\varphi$ is smooth on an open dense subset $U \subset H$ and $\varphi (U) \subset X_\reg$. Suppose that $Z$ is a closed subvariety of
 $H$ and $\codim_H (H\setminus U) > \dim Z$. Then for a general element $\alpha$ in a perfect
family $\cA$ of automorphisms of $H$ the morphism $\varphi|_{\alpha (Z)} : \alpha (Z) \to X$
is an injective immersion with the image in $X_\reg$. 
\eprop

\bproof  
By Theorem \ref{agga.t1} for a general element $\alpha$
in an algebraic family $\cA$ of automorphisms of $H$ the variety $\alpha (Z)$ does not meet $H\setminus U$.
By Proposition \ref{new.p1}(3)-(4) we can suppose that $\cA$ is perfect.
Theorem \ref{gp1.t1}(i)-(ii) and Proposition \ref{new.p1}(5)
imply now that for a general $\alpha \in \cA$ the morphism $\varphi|_{\alpha (Z)} : \alpha (Z) \to X$
is an injective immersion  which concludes the proof.
\eproof

\bthm\label{fle.t2} Let $Z$ be an affine algebraic variety and $X$ be a smooth quasi-affine flexible variety of dimension at least $\ED (Z)$.
Then $Z$ admits an injective immersion into $X$.
\ethm

\bproof    Let $U$, $H$, $H_i$ and $\Psi:  H\times X \to X\times X$  be as in Proposition \ref{pre.p2}, i.e., $H \simeq \A^t$ since each $H_i$ is a unipotent group.
Restricting $\Psi$    to $H\times x_0$ where $x_0$ is any point in $X$
we get a morphism $\varphi : \A^t \simeq H \to X$ which is smooth on $U$.                                                                                                                                                                                                     
By Holme's theorem we can treat $Z$ as a closed subvariety of $\A^t$. 
By Proposition \ref{pre.p2} we can suppose that $\codim_H H\setminus U >\dim Z$. 
 Since $\A^t$ is a smooth flexible variety
we get the desired conclusion by Proposition \ref{pre.p3}. 
\eproof

\brem\label{fle.r2} The author does not know if instead of ``injective immersion" one can use ``closed embedding" in Theorem \ref{fle.t2}.
Establishing properness is the bottleneck of the method known to the author. However, in the case of affine toric varieties
we managed to cope with this difficulty (see Theorem \ref{emb.t1} below).
\erem

\section{Affine toric varieits: preliminaries}
We suppose that readers are familiar  with toric varieties and all information about toric varieties which is used below can be found in the book of Cox, Little and Schrenck \cite{CLS}).
We fix  the following notations for  the rest of the paper.

- $N\simeq \Z^n$ - the standard lattice in $\R^n$;

- $M=\Hom_\Z (N, \Z)$ - the  lattice dual to $N$;

-  $<m,u>$  -  pairing of $m\in M$ or $M_\R=M\otimes_\Z \R$  
with $u\in N$ or $N_\R=N\otimes_\Z \R$;

- $\sigma$ - a rational convex polyhedral cone in $N_\R$;

- $\sigma^\vee$ - the dual cone of $\sigma$ in  $M_\R$;

- $\gamma^\bot$ - the set $\{ m \in M_\R| < m,\gamma>=0\}$  where $\gamma$ is any face of $\sigma$;

- $X_\sigma$ - the toric variety of $\sigma$, i.e., $X_\sigma$ is the spectrum of the group algebra of the semigroup $\sigma^\vee \cap M$;

- $\T=\Hom (M, \G_m)$ - the torus acting on $X_\sigma$; 

-- $\{ \rho_1, \ldots , \rho_r\}$ - the set of  extremal rays of $\sigma$  (by abusing notations we also denote by $\rho_i$ the ray generator, i.e., the 
primitive lattice vector on  the corresponding ray);

- $\sigma (k)$ - the set of $k$-dimensional faces of $\sigma$ (e.g., $\rho_i \in \sigma(1)$);

-  $O_i$ - the $\T$-orbit  in $X_\sigma$ corresponding to $\rho_i$ by the Orbit-Cone correspondence \cite[Theorem 3.2.6]{CLS};

- $D_i$ - the irreducible $\T$-invariant Weil divisor  in $X_\sigma$ containing $O_i$ as an open subset
(i.e., $D_i$ is the spectrum of the semigroup algebra of  $\tau_i=\rho_i^\bot \cap \sigma^\vee \cap M$);

- $H_i$ - the $\G_m$-subgroup of $\T$ corresponding to $\rho_i$, i.e., $H_i$ is a unique $\G_m$-subgroup of $\T$
that acts  trivially on $D_i$ and for $t \in H_i$ one has $t.\chi^m=t^{<m,\rho_i>}\chi^m$.

We would like to remind that since $\sigma^\vee  \cap M$ and $\tau_i=\rho_i^\bot \cap \sigma^\vee \cap M$ are  saturated affine semigroups the varieties $X_\sigma$ and $D_i$ are
always normal  (e.g., see \cite[Theorem 1.3.5]{CLS}).  Furthermore, we consider only the case when $X_\sigma$ {\em has no torus factors}
(or, equivalently, every invertible function on $X_\sigma$ is constant).

Let $X_\Sigma$ be the toric variety of a fan $\Sigma$ and $r$ be the cardinality of one-dimensional cones in $\Sigma$
(in  particular, if $\Sigma =\sigma$ then $r$ is the number of the ray generators $\rho_i$ of $\sigma$).
If torus factors are absent, then by \cite[Theorem 5.1.10]{CLS} there exists a subgroup $G$ of 
$\G_m^{r}$  which is a quasitorus (i.e., the direct product   of a finite subgroup and a subtorus of $\G_m^r$) 
and a closed subvariety $Z (\Sigma )$ of $\A^r$ such that  $Z (\Sigma )$ is invariant under the natural action of $G$ on $\A^r$ and
$X_\Sigma$ is isomorphic
to $(\A^{r} \setminus Z (\Sigma ))\sslash G$ (while $\G_m^{r}/G$ is isomorphic to the torus $\T$ acting on $X_\Sigma$). 
We are dealing with the situation when $\Sigma =\sigma$ and
in this case $Z(\sigma)$ is empty by construction (see the definition of   $Z(\sigma)$ on \cite[page 206]{CLS}). Thus, we have the quotient  morphism
\be\label{not.eq1} \pi : \A^{r}\to  \A^{r}\sslash G\simeq X_\sigma \ee which is $\G_m^r$-equivariant.
In connection with this formula we fix the following notations.

- $x_1, \ldots, x_r$ - a fixed coordinate system on $\A^r$;

- $\tT=\G_m^r$- the standard torus (with respect to the coordinate system) acting on $\A^r$;

- $\tD_i$ - the hyperplane in $\A^r$ given by $x_i=0$;

- $\tH_i$ - the $\G_m$-subgroup of $\tT$ acting trivially on $\tD_i$, i.e., this action is the flow of the semisimple vector field $x_i\frac{\p}{\p x_i}$;

- $U$ - the subset of $X_\sigma$ consisting of all points $u\in X_\sigma$ for which $\pi^{-1} (u)$ is a $G$-orbit
(in particular, this orbit is closed in $\A^r$);

- $U_0$ - the subset of $X_\sigma$ consisting of all points $u\in U$ for which the orbit $\pi^{-1} (u)$ has a trivial stabilizer.

Let us list some properties of the morphism $\pi : \A^r \to X_\sigma$ and the objects introduced before.

\blem\label{not.l1}
{\rm (i)} The morphism $\pi$ is an almost geometric quotient, i.e., $U$ is an open dense subset of $X_\sigma$ and,
consequently, general orbits of $G$ in $\A^r$ are closed and isomorphic to $G$.\footnote{Note that this implies that
$U_0$ is also open dense subset of $X_\sigma$ by \cite[Theorem 6.3]{PV}.} Furthermore, for every $u$ in $U$
the fiber $\pi^{-1} (u)$ is isomorphic as a homogeneous $G$-space to $G/F$ where $F$ is a finite subgroup of $G$.

{\rm (ii)}  The group $G\cap \tH_i$ is trivial. 

{\rm (iii)} Let $\theta\in \sigma (k)$ be regular, i.e., the set  $\{ \rho_{i_1}, \ldots, \rho_{i_k}\}$ of the generators
of the extremal rays of $\theta$
can be extended to a basis of $N$. Then $G$ meets the group $\tF$ generated by $\tH_{i_1}, \ldots, \tH_{i_k}$ at identity only.

 {\rm (iv)} The  homomorphism $\pi^* : \bk [X_\sigma] \to \bk [x_1, \ldots , x_r]$ induced by $\pi$ is determined by
the formula  \be\label{not.eq2} \pi^* (\chi^m)=\prod_{l=1}^r x_l^{<m,\rho_l>}.\ee

{\rm (v)} The image of $\tH_i$ in $\T=\frac{\tT}{G}$ coincides with $H_i$ and $\pi (\tD_i) =D_i$.

\elem

\bproof For the fist statement of (i) see  \cite[Theorem 5.1.10]{CLS}. This implies that  $\dim G=\dim A^r -\dim X_\sigma =r-n$.
For every $u \in U$ one has $\dim \pi^{-1} (u)=r-n$ since it cannot be less by
Chevalley's theorem  \cite[Chapter II, Exercise 3.22]{Har} and it cannot be larger since $\pi^{-1} (u)$ is a $G$-orbit. 
Being a homogeneous space $\pi^{-1}(u)$ is of the form $G/F$ where 
$F\subset G$ is the stabilizer of the orbit $\pi^{-1}(u)$ and the condition on the dimension implies that $F$ is finite. This concludes (i).

By \cite[Lemma 5.1.1]{CLS} we have 
{\be\label{not.eq3} G=\{\bar t=(t_1, \ldots, t_r) \in \G_m^r |\, \prod_{i=1}^r {t_i^{<\rho_i, m>}}=1 \text{ for all } m\in M\}.\ee}
Suppose that $G\cap \tH_i$ contains a finite subgroup of $d$-roots of unity. Assume that $d >1$ and let $\varepsilon$ be a 
primitive $d$-root of unity. Then Formula \eqref{not.eq3} implies that $\varepsilon^{<m,\rho_i>} =1$ for every $m\in M$,
i.e., $<m,\rho_i>$ is divisible by $d$. Hence, $\frac{\rho_i}{d} \in N$ contrary to the fact that $\rho_i$ is
a primitive vector in the lattice $N$. This yields (ii).

In (iii) let  $G\cap \tF$ contains a subgroup isomorphic to
the group of $d$-roots of unity. Then the similar argument implies that for a collection $\{ l_1, \ldots , l_s\}$ of integers
with greatest common divisor 1 the sum $\sum_{s=1}^k l_s \rho_{i_s}$ is divisible by $d$. However, if $d>1$, then this is contrary to the fact that
the set $\{ \rho_{i_1}, \ldots, \rho_{i_k}\}$ is extendable to a basis of $N$. Thus, we have (iii).

 For (iv) see \cite[page 209]{CLS}. 
Statement (v) follows from the explicit construction of $\pi$ in \cite[Proposition 5.1.9]{CLS} which implies, in particular (iv).
Vice versa, (v) can be also illustrated by
Formula \eqref{not.eq2}. Indeed, this formula implies that for every $\chi^m, \,  m \in \tau_i={\rho_i^\bot } \cap \sigma^\vee \cap M$
the function $\pi^*(\chi^m)$ is independent of $x_i$, i.e., it is fixed under the $\tH_i$-action. In particular, $\chi^m$ is fixed under the action of the image $H_i'$ of $\tH_i$ in $\T=\frac{\tT}{G}$
(which is isomorphic to $\tH_i$ by (ii)). Hence, the $H_i'$-action on $D_i$ is trivial. Since $H_i$ is a unique $\G_m$-subgroup of $\T$
with this property we see that $H_i'=H_i$. In particular, the divisor $\pi^{-1} (D_i)$ must be fixed under the $\tH_i$-action which implies
that $\pi^{-1} (D_i)=\tD_i$ and we are done.
\eproof

Consider an affine algebraic group $H$ acting on an affine variety $Z$, an affine variety $Y$
with a trivial $H$-action and an $H$-equivariant morphism $\rho : Z \to Y$.
Recall that under these assumptions
$Z$ together with the morphism $\rho : Z \to Y$  is called an $H$-torsor (or a principal $H$-bundle)
if for every $y \in Y$ there exists an \'etale morphism $\varphi_y : W_y \to Y$ such that
$y \in {\rm Im} \, \varphi_y$ and $W_y \times_YZ$ becomes a trivial principal $H$-bundle under the natural $H$-action.
If each $\varphi_y$ is injective, then $Z$ is called a locally trivial principal $H$-bundle.

\bprop\label{not.p1}

{\rm (1)} The morphism $\pi|_{\pi^{-1} (U_0)}: \pi^{-1} (U_0)\to U_0$
is a principal $G$-bundle (in particular,  $\pi$ is smooth over $U_0$). Furthermore,
if $G$ is irreducible, then this principal $G$-bundle
is locally trivial.

{\rm (2)}  Let $E$ be the subset of $\A^r$ consisting of all  $\bar x=(x_1, \ldots , x_r)\in \A^r$ such that at most one coordinate $x_i$ is equal to zero.
Then $\pi (E)\subset U_0$.

{\rm (3)} $U_0$ is the regular part of $X_\sigma$.

\eprop

\bproof   
If $u \in U_0$ and $w \in \pi^{-1} (u_0)$, then by the Luna \'etale slice theorem \cite{Lu} (see also \cite[Theorem 6.4]{PV}) there exists a smooth subvariety $V$ of $X_\sigma$
transversal to $\pi^{-1} (u)$ at $w$ such that $\pi|_V : V \to U_0$ is \'etale. This implies  that $U_0$ is contained in the regular part of $X_\sigma$
and that $\pi$ is smooth over $U_0$.
Furthermore, the natural $G$-action makes $V\times_{U_0}  \A^r$ a trivial $G$-bundle. Hence, $\pi^{-1} (U_0)\to U_0$ is a principal $G$-bundle.
Recall also that if $G$ is connected, then $G$ is a special group in the sense of Serre (see, \cite[Def. 2 and page 16]{Gro58}) and for any special
group $K$ every $K$-principal bundle is locally trivial \cite[Theorem 3]{Gro58}. Thus, we have (1).

For (2) and (3) we need to recall that by Orbit-Cone correspondence \cite[Theorem 3.2.6]{CLS} every $\theta \in \sigma (k)$ corresponds to a $\T$-orbit $O(\theta)\subset X_\sigma$
of dimension $n-k$ where $O(\theta)$ is the orbit of a so-called distinguished point. The description of this point \cite[page 116]{CLS}
implies that $O(\theta)$ consists of all points $u\in X_\sigma$ such that $\chi^m(u) \ne 0$ if and only if $m \in \theta^\bot \cap M$.
In particular, the ring  $\bk [R]$ of regular functions on the closure $R$ of $ {O(\theta)}$ in $X_\sigma$ can be viewed as the semigroup algebra of $ \theta^\bot \cap M$.
Let $\rho_{i_1}, \ldots, \rho_{i_k}$ be the extremal rays generating $\theta$ (i.e., $\theta^\bot=\rho_{i_1}^\bot\cap \ldots \cap \rho_{i_k}^\bot$)
and $F$ be the subgroup of $\T$ generated by $H_{i_1}, \ldots , H_{i_k}$.
Note that the natural inclusion $\bk [R]\hookrightarrow \bk[X_\sigma]$ makes $\bk[R]$ the subring of $F$-invariants and $R$ is given
in $X_\sigma$ by the ideal generated by  $\{ \chi^m | m \in (\sigma^\bot\setminus \theta^\bot) \cap M\}$. 
Hence, $R$ is the fixed point set of the $F$-action since for every $v\in X_\sigma \setminus R$ one can find $m \in (\sigma^\bot\setminus \theta^\bot)\cap M$
with $\chi^m(v) \ne 0$.
The difference between the points of $O(\theta)$ and $R\setminus O(\theta)$ is that for $w \in R\setminus O(\theta)$ there exists
$j\in \{1, \ldots, r\}\setminus \{i_1, \ldots, i_k\}$ such that $w$ is also fixed under the $H_{j}$-action
(because, $w$ is contained in a  $\T$-orbit of a smaller dimension corresponding to a cone in $\sigma$ containing $\theta$
and some $\rho_{j}$), whereas for any point in $O(\theta)$
such $j$ does not exist.

Let $\bar x\in \tT= E\setminus \bigcup_{i=1}^r\tD_i$.  Since $\pi$ is $\tT$-equivariant $G.\bar x$ is a general orbit, i.e., 
$\pi (\bar x) \in U_0$ by Lemma \ref{not.l1} (i). Since $E\cap \tD_{i_1}$ is a $\tT$-orbit dense in $\tD_{i_1}$ and $\pi (\tD_{i_1})=D_{i_1}$ by Lemma \ref{not.l1} (v)
we see that $\pi ( E\cap \tD_{i_1})=O_{i_1}$.
For every $\bar x \in E\cap \tD_{i_1}$ its $\tT$-orbit  $Q$ is naturally isomorphic to $G$ by Lemma \ref{not.l1} (ii) and $\pi (x) =u \in O_{i_1}$.
Note that if $Q$ is not closed, then its closure contain a point with some coordinates $x_{j}=0$ where $j\ne i_1$.
However, this implies that $u$ is a fixed point under both $H_{i_1}$-action and $H_{j}$-action contrary to the argument before.
Hence, $Q$ is closed. Furthermore, $Q$ is a unique closed $G$-orbit in $\pi^{-1} (u)$ by \cite[Theorem 5.0.7]{CLS}. 
If there exists another orbit $Q'$ in $\pi^{-1} (u)$, then the closure of $Q'$ contains $Q$ (e.g., see  \cite[Theorem 4.7 and Corollary]{PV}). However, 
this is impossible since $\dim Q' \leq \dim Q=\dim G$.  Hence, $\pi^{-1} (u)=Q$ and $u \in U_0$ which is (2).

Consider $u$ in the smooth part of $X_\sigma$. Then $u$ is contained in some $O(\theta)$ as before where $\theta$ must be regular
by \cite[Theorem 1.3.12 and Example 1.2.20]{CLS}.
Let $\tO(\theta)\subset \A^r$ be the $\tT$-orbit consisting of all points $\bar x$ whose zero coordinates
are exactly $x_{i_1}, \ldots, x_{i_k}$. Let $\theta'$ be a cone in $\sigma$ properly contained in $\theta$, i.e.,
$O(\theta)$ is contained in the closure of $O(\theta')$. Let us, say,
that $\theta'$ is generated by extremal rays  $\rho_{i_2}, \ldots, \rho_{i_k}$.
Then we can suppose by induction that such $\theta'$ is regular and that $\pi(\tO(\theta'))=O(\theta')\subset U_0$.
In particular,  $\pi^{-1} (u)$ belongs to the closure of $\tO(\theta')$  and $\pi^{-1} (u) \cap \tO(\theta')=\emptyset$.
This implies that every $\bar x \in \pi^{-1}(u)$ cannot have a nonzero coordinate
$x_{i_1}$. Hence, $\bar x $ must be contained in $\tO(\theta)$
(indeed, if $\bar x$ has a zero coordinate $x_j$ with $j \notin \{i_1, \ldots, i_k\}$, then $u$ is fixed under the $H_j$-action
contrary to the argument before). This implies that  the $\tT$-orbit $Q$ of $\bar x$ is closed since otherwise its closure
contains a point with an undesirable zero coordinate. By Lemma \ref{not.l1}(iii) $Q$ is naturally isomorphic to $G$
and arguing as before we see that $\pi^{-1}(u)=Q$.  Hence, $u \in U_0$ which yields (3) and concludes the proof. 
\eproof

\bcor\label{not.c1} Let $Y$ be a open subset of $X_\sigma$ such that ${\rm codim}_{X_\sigma} (X_\sigma\setminus Y)\geq 2$.
Then $\pi^{-1} (X_\sigma \setminus Y)$ has codimension at least 2 in $\A^r$.
\ecor
\bproof Note that $\pi^{-1} (Y) \subset \pi^{-1} (U\cap Y) \cup (\A^r \setminus E)$.
The definition of $U$ implies that $\pi^{-1} (U\setminus Y)$ has codimension at least 2 in $\A^r$ and the same is true for $A^r\setminus E$.
Hence, we have the desired conclusion.
\eproof

\bcor\label{not.c2} Let $C$ be a closed curve in $X_\sigma$ contained in  $U_0$.
Suppose that either

{\rm (1)}  $C$ is isomorphic to the affine line $\A^1$ or

{\rm (2)} $G$ is connected and $C$ is a smooth rational curve.

Then there exists a closed curve $\tC\subset \A^r$ such that $\pi|_\tC : \tC \to C$ is an isomorphism.

\ecor

\bproof  By Proposition \ref{not.p1} $\pi^{-1}(C)$ is a locally trivial principal $G$-bundle.
Statement (1) now follows from \cite[Theorem A.1]{FS} which states that  for each  affine algebraic group $F$
every principal $F$-bundle over the affine line  admits  a section. 

In (2) by Proposition \ref{not.p1}(1)  we can find an open cover $\{V_i\}$ of $C$ for which $\pi^{-1}(V_i)$ is naturally isomorphic to $V_i\times G$.  In particular, one has sections
$s_i : V_i\to \pi^{-1} (V_i)$ and $s_j|_{V_i\cap V_j}=g_{ij}s_i$ where $g_{ij} :V_i\cap V_j \to G$ is a morphism.
Since $G\simeq \G_m^{r-n}$ we see that $g_{ij}$ can be presented as a collection of $r-n$ sections of $\cO_{C}^*$ over $V_i\cap V_j$.
Hence, $H^1(C,G)$ is the direct sum of $r-n$ samples of $H^1(C, \cO_{C}^*)$. Since $C$ is a smooth rational curve 
we have $H^1(C, \cO_{C}^*)=\Pic C =0$ and, hence, $H^1(C,G)=0$. Thus, we can suppose that every pair of sections $s_i$ and $s_j$
agree on $V_i\cap V_j$.  Consequently, we have a global section of $\pi|_{\pi^{-1}(C)} : \pi^{-1}(C)\to C$ which yields the desired conclusion.
\eproof

\section{Embedding theorems for affine toric varieties} 

\bnota\label{emb.n1} 
In this section $X_\sigma$  is an affine toric variety without  torus factors.
In particular, $ X_\sigma \simeq \A^{r}\sslash G$ where $G\subset \tT =\G_m^r$ is a quasitorus acting naturally on $\A^r$.
We also denote by $ \pi : \A^{r}\to  X_\sigma $ the quotient morphism as in Formula \eqref{not.eq1}
with $U$ (resp. $U_0$) being the dense open subset of $X_\sigma$ consisting of all points $u \in X_\sigma$
for which $\pi^{-1} (u)$ is a $G$-orbit (resp. a $G$-orbit with a trivial stabilizer). That is, $U_0$ is the regular part of $X_\sigma$
by Proposition \ref{not.p1}.

\enota

\blem\label{emb.l1}
Let  
    $\bar 0$ be the origin in $\A^r$, $o=\pi (\bar 0)$ and $A=\bk[\A^r]^G$.          
Then one can choose a collection of monomials as generators of $A$ and the set $V$ of common zeros of this collection is 
contained in $\pi^{-1} (o)$. In particular,  $V\subset \pi^{-1} (X_\sigma \setminus U)$ unless $X_\sigma =U$
\footnote{Recall that if $X_\sigma =U$, then $\sigma$  is simplicial by \cite[Theorem 5.1.10]{CLS}.}.
\elem

\bproof  
Since the natural $\tT$-action respects monomials the same is true for the $G$-action.
Thus, any $G$-invariant polynomial is the sum of $G$-invariant monomials which yields the first claim.
Since $V\subset \A^r$ is closed and $G$-invariant $Z=\pi (V)$ is closed in $X_\sigma$ and
for every $z \in Z$ the only closed orbit in $\pi^{-1} (z)$ is contained also in $V$. 
Assume that $Z $ contains two distinct points $z_1$ and $z_2$ and $L_i$ is the closed orbit of $\pi^{-1} (z_i)$.
Note that the restriction of every polynomial from $A$ to $V$ is constant. Hence, elements of $A$ do not separate
$L_1$ and $L_2$ contrary to the fact that the regular functions on $X_\sigma$ separate $z_1$ and $z_2$. 
Thus, $Z$ is at most a singleton. Since $\bar 0\in V$ and $o \in Z$ we  see that $V\subset \pi^{-1} (o)$.
Let $\dim \pi^{-1} (z) > r-n$ for some $z \in X_\sigma\setminus o$. Since $\pi$ is $\tT$-equivarinat the same is true for all points in the $\T$-orbit $P$ of $z$ in $X_\sigma$.
The closure of $P$ contains a $\T$-orbit $Q$ of a smaller dimension and for every $w\in Q$ one has $\dim \pi^{-1} (w) >r-n$ by Chevalley's theorem.
Reducing the dimension of such $\T$-orbits further we see that $\dim \pi^{-1} (o) >r-n$, i.e., $o \notin U$. This concludes the proof.
\eproof

\bthm\label{emb.t1}  Let 
$X_\sigma$ be a normal affine toric variety without torus factors
and $l={\rm codim}_{ \A^r} \pi^{-1} (X_\sigma \setminus U_0)$.
Suppose that $Z$ is an affine algebraic variety such that $\ED (Z) \leq \dim X_\sigma$ and $\dim Z < l$.                      
Then there exists a closed embedding $\iota : Z \hookrightarrow X_\sigma$ such that $\iota (Z)$ is contained in the regular part $U_0$  of $X_\sigma$.
Furthermore, $l\geq 2$ and, in particular, for every affine curve $C$ with $\ED (C)\leq \dim X_\sigma$ there exists a closed embedding
of $C$  in  $X_\sigma$ with the image in $U_0$.
\ethm

\bproof   
 By Proposition  \ref{not.p1}(1)  the morphism $\pi|_{\pi^{-1}(U_0)}: \pi^{-1}(U_0)\to U_0$ is smooth
and $l \geq 2$ by Corollary \ref{not.c1}. By Holme's theorem $Z$ can be treated as a closed subvariety of $\A^r$. 
Proposition \ref{pre.p3} implies now that  for a general element $\alpha$ in a perfect family $\cA$ of automorphisms of $\A^r$
the morphism $\pi|_{\alpha (Z)} : \alpha (Z) \to U_0$ is an injective immersion.

Furthermore, consider the natural embedding $\A^r \hookrightarrow \PP^r$, $D=\PP^r \setminus \A^r\simeq  \PP^{r-1}$ and
$H=\GL_r(\bk)$. Then we have the natural $H$-action on $\PP^r$ such that $D$ is
invariant under it. By Proposition \ref{new.p1}(3) we can replace $\cA$
by the family $H\times \cA$. That is, for  a general $h$ in $H$ and a general $\alpha$ in $\cA$ the morphism
$\pi|_{h\circ \alpha (Z)} : h\circ \alpha (Z) \to U_0$ is still an injective immersion.

By Lemma \ref{emb.l1}  we can find generators $g_1, \ldots, g_s$ of $\bk [X_\sigma]$
such that the polynomials $f_i =g_i\circ \pi$
are monomials and the codimension (in $\A^r$)  of the variety given by $ f_1=\ldots = f_s=0$ is at least  $l$.
Note also  that $f_1, \ldots , f_s$ can be viewed as coordinate functions of $\pi : \A^r \to X_\sigma \subset \A^s$
and they can be extended to rational functions on $ \PP^r$. Since each $f_i$ is homogeneous with respect to the standard degree function
the intersection $R $ of the
indeterminacy sets  of these extensions is given by the common zeros of $f_1, \ldots,  f_s$ in $D$.
In particular, $R$ has codimension at least  $l$ in $D$. Let $P$ be the intersection of $D$ with the closure
of $h\circ \alpha (Z)$ in $\PP^r$, i.e., $\dim P \leq \dim Z -1< l-1$. 
Since the  restriction of the $H$-action to $D$ is transitive $P$ does not meet $R$  for
general $h \in H$ and $\alpha \in \cA$ by Theorem \ref{agga.t1}. 
Hence, $\pi|_{h\circ \alpha (Z)} : h\circ \alpha (Z) \to X_\sigma$ is a proper morphism
by \cite[Corollary 5.4]{Ka20}.  Consequently, it is a closed embedding which concludes the proof.
\eproof

In particular, we have the following fact which is also a trivial consequence of Corollary \ref{fle.c1}.

\bcor\label{emb.c1} Let $X_\sigma$ be an affine simplicial toric variety. Let $Z$ be an affine algebraic variety such that $\ED (Z)\leq \dim X_\sigma$ 
and $\dim Z$ is less than
the codimension of the singularities of $X_\sigma$ in $X_\sigma$.  Then there is a closed embedding of $Z$ into
$X_\sigma$ with the image in $U_0$.
\ecor

\section{Locally nilpotent vector fields on affine toric varieties}\label{lnd.s1}

We use a combinatorial description of locally nilpotent vector fields on $X_\sigma$ given by Liendo in his paper \cite{Li}
in which he rediscovered  Demazure roots  \cite[Section 3.1]{De}  (see also \cite[Definition 4.2]{AKuZ19}).
Recall that a Demazure root associated with some $\rho_i$ is any element $e$ of $M$ such that $<e,\rho_i>=-1$ and  $<e,\rho_j>$ 
is nonnegative for every $j\ne i$.  The vector field on $X_\sigma$ defined by
\be\label{lnd.eq1} \p_{\rho_i,e} (\chi^m) =<m, \rho_i>\chi^{m+e} \ee
is locally nilpotent and up to a constant factor every homogeneous locally nilpotent vector field is of this form.
For a  Demazure root $e \in M$ associated with $\rho_i$
one has $\tilde e: = \pi^* (\chi^e) =(\tilde e_1, \ldots , \tilde e_r)$ where by Formula \eqref{not.eq2}  the  $i$-th coordinate $\tilde e_i$
is equal to $-1$. Let  $\tilde e' =(\tilde e_1', \ldots , \tilde e_r')$ where the  $i$-th coordinate $\tilde e_i'$
is equal to zero, whereas $\tilde e_l'=\tilde e_l$ for $l\ne i$.
Formulas \eqref{not.eq2} and \eqref{lnd.eq1} imply now the following fact which was first discovered in \cite{AKuZ19}).

\blem\label{lnd.l1} The polynomial $\pi^* (\p_{\rho_i,e} (\chi^m))$ coincides with $\tilde \p_{\rho_i,e} (\pi^*(\chi^m))$
where  the locally nilpotent vector field $\tilde \p_{\rho_i,e}$ on $\A^r$  is given by $\bar x^{\tilde e'}\frac{\p}{\p x_i}$
with $\bar x^{\tilde e'}=\prod_{l=1}^r x_l^{\tilde e_l'}$, i.e., the flow of $\tilde \p_{\rho_i,e}$ is given by
\be\label{lnd.eq2} \bar x =(x_1, \ldots, x_r) \mapsto (x_1, \ldots, x_{i-1}, x_i+t\bar x^{\tilde e'}, x_{i+1}, \ldots , x_r) \ee
where $t$ is the time parameter.
\elem

The algebra $\bk[D_i]$ of regular functions on $D_i$ can be viewed as the semigroup algebra of
$\tau_i=\rho_i^\bot\cap \sigma^\vee \cap M$. 
Note that $\bk[D_i]$ is the kernel of $\p_{\rho_i,e}$ viewed as a derivation on $\bk[X_\sigma]$. 
The natural embedding $\bk[D_i] \hookrightarrow \bk [X_\sigma]$ yields a dominant $\T$-equivariant morphism $\kappa_i : X_\sigma \to D_i$
that is the categorical quotient of the $\G_a$-action associated with $\p_{\rho_i,e}$. Note that it is also the categorical quotient of the natural $H_i$-action
on $X_\sigma$ since $\bk[D_i]$ is the subring of $H_i$-invariants of $\bk [X_\sigma]$.
Furthermore, as we mentioned before in the proof of Proposition \ref{not.p1} $D_i$ is the fixed point set of the $H_i$-action  on $X_\sigma$.

\bnota\label{lnd.n1} Similarly, consider a cone $\theta \in \sigma (2)$ containing two extremal rays $\rho_i$ and $\rho_j$
and the subgroup  $H_{ij}$ of $\T$ generated by $H_i$ and $H_j$. The dual cone of $\theta$ meets $M$ along
$\tau_i\cap\tau_j$. The semigroup algebra of $\tau_i\cap\tau_j$ can be viewed as the algebra of regular functions
on $D_{ij}=D_i\cap D_j$. As before, one can see that $D_{ij}=X_\sigma\sslash H_{ij}$ and $D_{ij}$ is the fixed point set
of the $H_{ij}$-action on $X_\sigma$.  Since $\bk [D_{ij}]$ has no zero-divisors and its transcendence degree is $n-2$
one can see that $D_{ij}$ is an irreducible $\T$-invariant Weil divisor in $D_i$.  In particular, $D_{ij}$ contains a dense
$\T$-orbit $O(\theta)$ (which is associated with $\theta$ via the Orbit-Cone correspondence).
\enota

\blem\label{lnd.l2} Let Notation \ref{lnd.n1} hold and $\theta$ be regular. Then $\kappa_i$ is smooth over  $O(\theta)$
and $\kappa_i^{-1}(u)$ is isomorphic to $\A^1$ for every $u \in O(\theta)$.
\elem

\bproof
Since $\theta$ is regular 
$O(\theta)$ is contained in the regular part $U_0$ of $X_{\sigma}$ by \cite[Theorem 1.3.12 and Example 1.2.20]{CLS}.
Let $u \in O(\theta)$. Then $T_uX_\sigma$ is equipped with the induced linear $H_{ij}$-action.
By the  Luna slice \'etale theorem for smooth points (e.g., see \cite[Theorem 6.4]{PV}) 
there exists
an $H_{ij}$-equivariant \'etale morphism $\varphi : Y \to T_uX_\sigma$ from a dense open $H_{ij}$-invariant subset $Y$ of $X_\sigma$ containing $u$.
Hence, since the map $T_uX_\sigma \to T_u X_\sigma\sslash H_i$ is smooth so is $\kappa_i|_Y : Y\to Y\sslash H_i$.
For every point $w \in \kappa_i^{-1}(u)$ the closure of its $H_i$-orbit contains the fixed point $u$ (e.g., \cite[Theorem 4.7 and Corollary]{PV}), i.e.,
 this orbit is contained in $Y$. Hence, $\kappa_i^{-1}(u)\subset Y$ (and, consequently, $\kappa_i^{-1}(u)$ is isomorphic via $\varphi$
 to an affine line through the origin in $T_uX_\sigma$). This yields the desired conclusion.
\eproof

\blem\label{lnd.l3} Suppose that for  some $\rho_i$ every $\theta\in \sigma(2)$ containing $\rho_i$ is regular.
Then there exists an open subset $V_i$ of $D_i\cap U_0$ such that $\codim_{D_i}D_i\setminus V_i\geq 2$ and for every $v \in V_i$
one can find a locally nilpotent vector field $\delta$ of the form $g\p_{\rho_i,e}$ where $g \in  {\bk (D_i)\subset \bk (X_\sigma)}$
which does not vanish on $\kappa_i^{-1} (v)$.
\elem

\bproof  
Let $\rho_{j_1}, \ldots, \rho_{j_k}$ be the collection of all extremal rays distinct from $\rho_i$ such that
for every $s=1, \ldots, k$ there exists $\theta \in \sigma(2)$ containing $\rho_{j_s}$ and $\rho_i$. 
Formula \eqref{lnd.eq1} implies that $\p_{\rho_i,e}$ does not vanish over $O_i=D_i \setminus \bigcup_{s=1}^kD_{j_s}$.
Thus, we have to consider the case when $v$ is a general point of some $D_{j_s}$.
Choose a rational function $f_s$ on $D_i$ with poles on $D_i \cap D_{j_s}$ only such that these poles are simple
at general points of $D_i \cap D_{j_s}$. Let $l_s$ be the zero multiplicity of $\p_{\rho_i,e}$ at general points of $D_{j_s}$.
Then the vector field $\delta=f_s^{l_s}\p_{\rho_i,e}$
is regular, locally nilpotent, tangent to the fibers of $\kappa_i$ and it does not vanish at general points of $D_{j_s}$.  
By Lemma \ref{lnd.l2}
$\kappa_i^{-1}(v)$ is isomorphic to the affine line and since $v \in D_i\cap D_{j_s}$ one has
$\kappa_i^{-1}(v)\subset D_{j_s}$.  Thus,   $\delta|_{\kappa_i^{-1}(v)}$
does not vanish since it is tangent to the  line $ \kappa_i^{-1}(v)$ and nonzero at a general point of  $ \kappa_i^{-1}(v)$. This yields the desired conclusion.
\eproof

\bprop\label{lnd.p1} Let every $\theta\in \sigma (2)$ containing $\rho_i$ be regular and let $V_i\subset D_i$ be as in Lemma \ref{lnd.l3}. 
 Let $Z$ be a closed subvariety
of $D_i\subset X_\sigma$ which is contained in $V_i$. Let $s : Z \to \A^r$ be a section of $\pi : \A^r\to X_\sigma$ over $Z$ for 
which $\tZ=s(Z)$ is closed in $\A^r$. Then one can find a locally nilpotent vector field $\delta$ equivalent to $\p_{\rho_i,e}$
\footnote{Two locally nilpotent derivations are equivalent  if they have the same kernels.} 
and such that $\delta$ does not vanish on $\kappa_i^{-1} (Z)$.
\eprop

\bproof  Lemma \ref{lnd.l3} implies that for every $ { u \in} Z$ one can find a locally nilpotent vector field $\delta_z$ 
of the form $g_z\delta_{\rho_i,e}, \, g_z \in  {\bk (D_i)\subset \bk(X_\sigma)}$ which does not vanish on $\kappa_i^{-1}(u)$. 
Recall that $g_z$ as an element of $ {\bk (X_\sigma)}$ is invariant under the $H_i$-action.
Hence, by Lemma \ref{lnd.l1}
$\delta_z=\pi_*(\tilde \delta_z)$ where $\tilde \delta_z$ is of the form $\tilde f_z \frac{\p}{\p x_i}$ with $\tilde f_z$ being
a polynomial independent of $x_i$ since $x_i$ is not invariant under the $\tH_i$-action.  
By assumption $\tilde \delta_z$ and, therefore, $\tilde f_z$ do not vanish
at $\pi^{-1}  {(z)}\cap \tZ$. Hence, by the Nullstellensatz one can find polynomials $\tilde h_z$
such that only finite number of them are nonzero and $\sum_z \tilde h_z\tilde f_z|_\tZ=1$.
By the assumption every regular function on $\tZ$ is a lift of a regular function on $Z$ which extends to an element $\bk [D_i]=\Ker \p_{\rho_i,e}\subset \bk [X_\sigma]$.
In particular, one can suppose that  $\tilde f_z=\pi^*f_z$ and $\tilde h_z=\pi^*h_z$ where $f_z,h_z \in \bk[D_i]$.
Hence, $\delta=\sum_zh_z\delta_z$ (resp. $\tilde \delta =\sum_z\tilde h_z\tilde \delta_z$) is a locally nilpotent non-vanishing  vector field on $\kappa_i^{-1} (Z)$
(resp. $\pi^{-1} (\kappa_i^{-1} (Z))$) which yields the desired conclusion.
\eproof

\bcor\label{lnd.c1} Let the assumptions of Proposition \ref{lnd.p1} hold and $Z$ (and, therefore, $\tZ$) be isomorphic to the affine line $\A^1$ equipped with
a coordinate $t$.  Let $\tilde \delta$ be the locally nilpotent vector field on $\A^{r}$ as in the proof of Proposition \ref{lnd.p1} (i.e.,  $\pi_*(\tilde \delta)=\delta$).
Then for every polynomial $h(t)$  there exists a function $g\in \bk [D_i]\subset \bk [X_\sigma]$
such for the flow $\tilde \beta_h^i$ of the locally nilpotent vector field  $\pi^*(g)\tilde \delta$ a time 1 \footnote{That is,
$\tilde \beta_h^i: \A^r\to \A^r$ is defined by $(\tilde \beta_h^i)_* (\lambda)=\exp (\pi^*(g)\tilde \delta) (\lambda) \, \forall \lambda \in \A^{[r]}$.}
 one has $x_i\circ \tilde \beta_h^i (t)=h(t), \, t \in \tZ$.
\ecor

\bproof Since $\delta$ is equivalent to $\p_{\rho_i,e}$ and does not vanish on $Z$ one can suppose  (by Lemma \ref{lnd.l1}) that
the restriction of $\tilde \delta$ to $\tZ$ coincides with $\frac{\p}{\p x_i}$.
Let $\check g(t)=h(t) -x_i(t)$.  Note that $\check g(t)$ (as a function on $Z$)   admits an extension to a function $g \in \bk [D_i]=\Ker \delta$.  
This extension yields the desired function.
\eproof

\section{Affine simplicial toric varieties}

Recall that an affine toric variety $X_\sigma$ is simplicial if every face of $\sigma$ is a simplex, i.e., $n=r$ and
$G$ is a finite group. This implies, in particular, that for every $j\ne l$
the extremal rays $\rho_l$ and $\rho_j$ are contained in some $\theta \in \sigma(2)$ and
 $D_{lj}=D_l \cap D_j$ is  always a Weil divisor in $D_l$.

\blem\label{sim.l1}  Let $X_\sigma$ a simplicial toric variety of dimension at least 4 which is
smooth in codimension 2 and $C$ be a smooth polynomial curve in the regular part of $X_\sigma$.
Let $V_i$ be as in Lemma \ref{lnd.l3}, 
$V=\bigcap_{i=1}^r\kappa_i^{-1}(V_i)$, $W_{l}=\kappa_l^{-1} (V_l) \setminus \bigcap_{j\ne l}  {\kappa_j^{-1}(V_j)}$ 
and $W_{l}'=\kappa_l(W_l)$.
For every $\theta \in \sigma (2)$  containing extremal rays $\rho_l$ and $\rho_j$ let
$\psi_{\theta} : X_\sigma \to D_{lj}$ be the morphism induced by the homomorphism of the semigroup algebras
associated with the natural embedding $\tau_l\cap \tau_j \hookrightarrow \sigma$.
Then replacing $C$ with its image under an automorphism of $X_\sigma$ one can suppose that

{\rm (i)} $C$ is contained in $V$;

{\rm (ii)}  $C_l=\kappa_l (C)$  meets $W_l'$ at a finite set for every $l=1, \ldots , r$;

{\rm (iii)} $\kappa_l|_C : C \to D_l$ is a closed embedding for every $l=1, \ldots, r$;

{\rm (iv)} $\psi_\theta : C \to  {\psi_\theta (C)}$ is a birational morphism
for every $\theta \in \sigma(2)$ containing $\rho_l$ and $\rho_j$.

\elem

\bproof Since $D_{lj}$ is a divisor in $D_l$ and  $D_l\setminus V_l$ has codimesion at least 2 in $D_l$ we see that
$V_l$ contains an open subset of $D_{lj}$ and $\kappa_l^{-1}(V_l)$ contains an open part of $D_j$.
Hence, $X_\sigma\setminus \kappa_l^{-1} (V_l)$ does not contain Weil divisors in $X_\sigma$, i.e., it is of codimension at least 2.
Consequently, ${\rm codim}_{X_\sigma} X_\sigma \setminus V \geq 2$.  Recall that $X_\sigma$ is flexible by \cite{AKuZ}
and, therefore, $U_0$ and $V$ are flexible by \cite[Theorem 2.6]{FKZ}. By Theorem \ref{agga.t1} for a general $\alpha$ in any perfect family $\cA$ of
automorphisms of $U_0$ (which are extendable to automorphisms of $X_\sigma$ by the Hartogs' theorem)
$\alpha (C)$ is contained in $V$ and $\alpha (C)$  meets every $\kappa_l^{-1}(W_l')$ at a finite set which yields (i) and (ii). 
Lemma \ref{lnd.l2}  implies that every $\kappa_l$ is smooth over $\kappa_l(V)$. Thus, by Theorem \ref{gp1.t1} and Proposition \ref{new.p1}(5) 
for a general $\alpha\in \cA$ each morphism $\kappa_l: \alpha(C) \to D_l$ is a closed embedding 
and each morphism $\psi_\theta: \alpha(C) \to {\psi_\theta (\alpha (C))}$ is birational which yields (iii)-(iv) and the desired conclusion.
\eproof

\blem\label{sim.l2}  Let the assumptions of Lemma \ref{sim.l1} hold and $C$ satisfy conditions (i)-(iv).
Suppose that $i, \delta$, $h$ and $g$  are as in Corollary \ref{lnd.c1} and $\beta_h^i$ is  flow of the locally nilpotent vector field  $g \delta$ a time 1.
Suppose further that $h(t)=ct+d$ where $c$ and $d$ are general constants. Then the curve  $\beta_h^i (C)$
also satisfies conditions (i)-(iv).
\elem

\bproof  Let $\tV =\pi^{-1}(V)$. By Corollary \ref{not.c2} there exists a curve $\tC\subset \tV$ such $\pi|_\tC : \tC \to C$ is an isomorphism.
Note that $\pi \circ \tilde \beta_h^i=\beta_h^i \circ \pi$ where $\tilde \beta_h^i$ is as in Corollary \ref{lnd.c1}.
Hence, besides conditions (i) and (ii) for  $\beta_h^i (C)$  it suffices to prove that

{\rm (iii$'$)}  for $\tilde \kappa_l= \kappa_l\circ \pi$ the morphism $\tilde \kappa_l|_{\tilde \beta_h^i(\tC)} : \tilde \beta_h^i(\tC) \to D_l$ is a closed embedding for every $l=1, \ldots, r$;

{\rm (iv$'$)} for $\tilde \psi_\theta =\psi_\theta \circ \pi$ the morphism  
$\tilde \psi_\theta|_{\tilde \beta_h^i(\tC)} : \tilde \beta_h^i(\tC) \to { \tilde \psi_\theta ({\tilde \beta_h^i(\tC)}) }$ is  birational
for every $\theta \in \sigma(2)$ containing $\rho_l$ and $\rho_j$.

Let us start with  (iv$'$).  One can choose  coordinate functions of $\tilde \psi_\theta$ in the form $\pi^* (\chi^m)$ where $m \in \tau_l\cap  \tau_j$.
By Formula\eqref{not.eq2} $\pi^* (\chi^m)$ is of the form $x_i^{k_m}y_m$ where $y_m$ is a monomial independent of $x_i$.
Condition (iv) implies that there exist  $m',m'' \in \tau_l\cap  \tau_j$ such that for $t\in \A^1\simeq \tC$ the functions
$ x_i^{k_{m'}}(t)y_{m'}(t)$ and $x_i^{k_{m''}}(t)y_{m''}(t)$ are not proportional and, in particular, $\frac{y_{m'}(t)}{y_{m''}(t)}$ is a nonzero rational function.
Hence, for general $c$ and $d$ the functions $ (ct+d)^{k_{m'}}y_{m'}(t)$ and $(ct+d)^{k_{m''}}(t)y_{m''}(t)$ are not proportional
and Corollary \ref{lnd.c1} implies that the morphism 
$\tilde \psi_\theta|_{\tilde \beta_h^i(\tC)} : \tilde \beta_h^i(\tC) \to { \tilde \psi_\theta ({\tilde \beta_h^i(\tC)}) }$ is  birational which is (iv$'$).

Let $S_{lj}$ be a finite subset of $\tC$ for which $\tilde \psi_\theta|_{\tilde \beta_h^i(\tC\setminus S_{lj})} : \tilde \beta_h^i(\tC\setminus S_{lj}) \to D_{lj}$
is an embedding. Then condition (iii) implies that for every $t_0\in S_{lj}$ there exists  $m\in { \tau_l}$ such that
$\frac{\dd}{\dd t} x_i^{k_m}y_m|_{t=t_0}$ is nonzero. This implies that either $y_m(t_0)\ne 0$ or $\frac{\dd}{\dd t}y_m|_{t=t_0}\ne 0$.
Consequently, $\frac{\dd}{\dd t} (ct+d)^{k_m}y_m|_{t=t_0}$ is nonzero for general $c$ and $d$. Hence, we can suppose that
$\tilde \kappa_l|_{\tilde \beta_h^i(\tC)} : \tilde \beta_h^i(\tC) \to D_l$ is an immersion. Furthermore, for every $t_0\ne t_1 \in S_{lj}$
there exists  $m\in  { \tau_l}$ such that $x_i^{k_m}(t_0)y_m(t_0) \ne x_i^{k_m}(t_1)y_m(t_1)$. Again for general $c$ and $d$
this implies that  $(ct_0+d)^{k_m}y_m(t_0) \ne (ct_1+d)^{k_m}y_m(t_1)$. Hence, $\tilde \kappa_l|_{\tilde \beta_h^i(\tC)} : \tilde \beta_h^i(\tC) \to D_l$ is a
closed embedding which is (iii$'$).  

By (ii) for every $l$  the curve $C_l$ meets $W_l'$ at a finite set $Q_l'$ or, equivalently,
for a general point $z \in C$ one has $\kappa_l(z)\notin W_l'$. The same remains true for a general point $\beta_h^i (z)$  of the curve
$\beta_h^i(C)$. Indeed, it suffices to show that it is true for some point of $\beta_h^i (C)$. Let $t\in \tC\simeq \A^1$ be the preimage of $z$ in $\tC$.
Note that $\beta_h^i(C)$ meets $C$ at the points where $x_i(t)=ct+d$. Since $c$ and $d$ are general the solutions of the latter equation
yield general points of $C$ and, hence,  condition (ii) for the curve  $\beta_h^i (C)$.

Recall that $\beta_h^i$ is the flow at time 1 of a locally nilpotent vector field  $g\delta$ as in Corollary \ref{lnd.c1} which is equivalent to a vector field $\p_{\rho_i,e}$
and, therefore,  { tangent} to the fibers of $\kappa_i$. 
In particular, $\beta_h^i(C)$ is contained in $\kappa_i^{-1} (C_i)$ where $C_i=\kappa_i (C)$.
By construction $\kappa_i^{-1}( C_i \setminus Q_i')\subset V$.
By Lemma \ref{lnd.l2} every fiber $L$ of $\kappa_l|_{\kappa_l^{-1}(Q_i')}: \kappa^{-1} (Q_i') \to Q_i'$ is isomorphic to the affine line and
such $L$ meets  $V$ since $C$ does at some point $z_0\in C$ (where this $z_0$ is unique since $\kappa_l : C \to D_l$ is a closed embedding).  
Let $t$ be a coordinate on $\tC\simeq C\simeq \A^1$.
Recall that by construction in Corollary \ref{lnd.c1} $g\subset \bk [D_i]$ is an extension of the function $ct+d-x_i(t)$.
Hence, $g\delta =\delta_1+ d\delta$ where the locally nilpotent vector field $\delta_1$ commutes with $\delta$. In particular, 
the flow of $g\delta$ at time 1
 { is} the composition of the flows of $\delta_1$ at time 1 and $\delta$ at time $d$. Since by Proposition \ref{lnd.p1} $\delta$ does not
vanish on $L$ and $d$ is general we see that $\beta_h^i(z_0)$ is a general point of $L$ and, therefore, it belongs to $V$.
This yields condition (i) for $\beta_h^i(C)$ and the desired conclusion.
\eproof

\bthm\label{sim.t1}  Let $X_\sigma$ be an affine simplicial toric variety  of dimension at least 4 such that $X_\sigma$ is
smooth in codimension 2. Suppose that $\T_0$ is an algebraic torus and $Y=\T_0\times X_\sigma$.
Let  $\varphi : C \to C'$ be an isomorphism of two smooth polynomial curves contained
in the regular part of $Y$. Then $\varphi$ extends to an automorphism of $Y$. 
\ethm

\bproof  
Let $\nu : Y \to \T_0$ be the natural projection. Since  $C$ and $C'$ are polynomial curves their images $z=\nu (C)$
and $z'=\nu (C')$ are singletons.
Consider an automorphism $\alpha_0$ of $\T_0$ sending $z'$ to $z$ and its natural lift to
an automorphism $\alpha$ of $Y$. Replacing $C'$ by $\alpha (C')$ and $\varphi$ by  $\varphi\circ \alpha^{-1}$ we
can suppose that $z=z'$. Hence, $C, C'\subset \nu^{-1}(z) \simeq X_\sigma$ and it suffices to consider the case of $Y=X_\sigma$ only.

Suppose that $\pi : \A^r \to X_\sigma$ is as in Formula \eqref{not.eq1}. By Corollary \ref{not.c2} one can find a curve        
$\tC$ (resp. $\tC'$)  in $\A^r$ such that $\pi|_{\tC}: \tC \to C$ (resp. $\pi|_{\tC'}: \tC' \to C'$) is an isomorphism.
Let $t'$ be a coordinate on $C'\simeq \tC'$ and $t=\varphi^*(t')$ be the coordinate on $C\simeq \tC$.
Applying consequently automorphisms $\beta_h^i$ as in Lemma \ref{sim.l2} with $i$ running over $\{1, \ldots, r\}$
one can suppose that $\tC$ is a curve such that $x_i(t)=c_it+d_i$ for every $i$ where $(c_1,d_1, \ldots, c_r,d_r)$ is a general point in $\A^{2r}$.
Similarly, one can suppose that $\tC'$ is a curve such that $x_i(t')=c_i't'+d_i'$  for every $i$ where $(c_1',d_1', \ldots, c_r',d_r')$ is a general point in $\A^{2r}$.
Choosing these two general points equal we get the desired conclusion.
\eproof

We need to remind the following \cite[Definition 8.2]{KaUd}.

\bdefi\label{sim.d1} 
Let  $C_1$ and $C_2$ be smooth  curves
in a smooth quasi-affine variety $Y$ with defining ideals $I_1$ and $I_2$ in $\bk [Y]$. 
We suppose also that $C_1$ and $C_2$ are closed in an affine variety containing $Y$.
Let $Y$ possess  a volume form $\omega$ (i.e., $\omega$ is a nonvanishing section of the canonical bundle on $Y$) and let
each conormal bundle $\frac{I_j}{I_j^2}$ of $C_j$ in $Y$ be trivial. 
By \cite[Lemma 6.3]{KaUd} there is a neighborhood $W_j$ of $C_j$ in $Y$ in which  $C_j$ is a strict complete intersection 
given by $u_{1,j}=\ldots =u_{n-1,j}=0$ where  $u_{1,j}, \ldots , u_{n-1,j}\in I_j$ and $n=\dim Y$. That is, for  $A_j =\frac{\bk[Y]}{I_j}$ we have the graded algebra
$\frac{\bk[Y]}{I_j^k}\simeq A_j\oplus \bigoplus_{l=1}^{k-1}\frac{I_j^l}{I^{l+1}_j}$
which can be viewed as the algebra of polynomials in $u_{1,j}, \ldots , u_{n-1,j}$ over $A_j$ of degree at most $k-1$.
Consider an isomorphism
 $\varphi : \frac{\bk[Y]}{I_1^k} \to \frac{\bk[Y]}{I_2^k}$  of these algebras for a natural $k$. 
Up to the induced isomorphism $A_1\simeq  A_2$ this isomorphism 
$\varphi$ is determined by its values $\varphi (u_{i,1}), \, i=1, \ldots, n-1$.
These values can be viewed as polynomials in $u_{1,2}, \ldots , u_{n-1,2}$ over $A_2$, i.e., one has the matrix $\left[\frac{\p \varphi (u_{l,1})}{\p u_{s,2}}\right]_{l,s=1}^{n-1}$.
 Since  the normal bundle $N_Y C_j$ is trivial, the existence of $\omega$ implies the existence  of a volume form on $C_j$.
Fix  volume forms $\omega_j$ on $C_j$ such that $\tilde \varphi^* \omega_1=\omega_2$ where 
the isomorphism $\tilde \varphi : C_2\to C_1$ is induced by 
$\varphi$. Choose a section $\pr_j : TY|_{C_j}\to TC_j$ of the canonical inclusion $TC_j\hookrightarrow TY|_{C_j}$ 
and consider the section $\tilde \omega_j=\omega_j\circ \pr_j $ of the dual bundle $(TY|_{C_j})^\vee$ of $TY|_{C_j}$.
Then one can require that
 $\omega|_{C_j}$ coincides with $\tilde \omega_j \wedge \dd u_{1,j} \wedge \ldots \wedge \dd u_{n-1,j}$.
Under this requirement  the determinant of  $\left[\frac{\p \varphi (u_{l,1})}{\p u_{s,2}}\right]_{l,s=1}^{n-1}$ 
is  well-defined modulo $I_2^{k-1}$ 
(i.e., it is independent of the choice of coordinates $u_{1,j}, \ldots , u_{n-1,j}$).
Hence, we say that $\varphi$ has Jacobian $a\in \bk \setminus \{ 0\}$ if the determinant of  $\left[\frac{\p \varphi (u_{l,1})}{\p u_{s,2}}\right]_{l,s=1}^{n-1}$ is equal to $a$ modulo $I_2^{k-1}$.
\edefi

Note that Definition \ref{sim.d1} is applicable in the case when $C_1$ and $C_2$ are smooth polynomial curves in a simplicial toric variety $X_\sigma$
contained in its regular part. Indeed,  $U_0$ as the regular part of $X_\sigma$ is flexible \cite{AKuZ}.
Recall that $\A^r$ admits a volume form invariant under the natural $\SL_r(\bk)$-action.
Hence, we can push this volume form down to  a volume form $\omega$ on $U_0$  since $\pi|_{\pi^{-1} (U_0)}: \pi^{-1} (U_0)\to U_0$ is an unramified covering by Proposition \ref{not.p1}.
Furthermore, since the normal bundles of smooth polynomial curves are always trivial,  we are under the assumptions of   Definition \ref{sim.d1}.

\bcor\label{sim.c1}   Let  $\varphi : \cC \to \cC'$ be an isomorphism of $k$th infinitesimal neighborhoods of two smooth polynomial curve  $C$ and $C'$ contained
in the regular part $U_0$ of an affine simplicial toric variety $X_\sigma$ of dimension at least 4 which is smooth in codimension 2.
Suppose that the Jacobian of $\varphi$ is a nonzero constant $a$.
Then $\varphi$ extends to an automorphism $\Phi$ of $X_\sigma$.
\ecor

\bproof  Recall that by \cite[Lemma 6.2]{KaUd}  every automorphism $\alpha$ of $U_0$ has a constant Jacobian where the Jacobian is computed as $\frac{\alpha^*(\omega)}{\omega}$
and if $\alpha$ is a composition of elements of flows of locally nilpotent vector fields, then its Jacobian is 1.
Let $\psi : C'\to C$ be an isomorphism.
By Theorem \ref{sim.t1} $\psi$ extends to an automorphism $\Psi$ of $X_\sigma$ which in turn induces an automorphism $\cC' \to \cC$ 
also denoted by $\psi$.  By construction $\Psi$ is a composition of elements of flows of locally nilpotent vector fields. Hence, its Jacobian is 1.
Taking a composition of $\varphi$ with the action of an appropriate element of $\T$  and replacing 
$C'$ with its image under this action we can suppose that the Jacobian of $\varphi$ is 
also 1 nodulo $(I')^{k-1}$ in the sense of Definition \ref{sim.d1} where $I'$ (resp. $I$) is the defining ideal of $C'$  (resp. $C$) in $\bk [X_\sigma]$.
Then the automorphism $ \lambda:=\psi \circ \varphi : \cC \to \cC$ has Jacobian 1 modulo $I^{k-1}$. By \cite[Theorem 6.6]{KaUd}
$\lambda$ extends to an automorphism $\Lambda$ of $X_\sigma$.
It remains to note that $ \Psi^{-1} \circ \Lambda$ is the desired extension of $\varphi$ and we are done.
\eproof

\brem\label{sim.r1} Let $\cN$ be the smallest saturated
set of locally nilpotent vector fields on $X_\sigma$ that contains all vector field of the form $\p_{\rho_i,e}$ as in Formula \eqref{lnd.eq1}.
Consider the subgroup $G\subset \SAut (X_\sigma)$ generated by $\cN$.
Assume that in Theorem \ref{sim.t1} $Y=X_\sigma$.  Then it follows from the proof that an automorphism
extending $\varphi$ can be chosen in $G$. Furthermore, one can check that $V$ as in Lemma \ref{sim.l1}  is $G$-flexible
(recall that we can suppose that such $V$ contains $C$ and $C'$). 
Hence, if the  Jacobian of $\varphi$ from Corollary \ref{sim.c1} is 1,
then \cite[Theorem 6.6]{KaUd} implies that $\Phi$ can be also chosen in $G$.         

\erem


\begin{thebibliography}{KaMi} 

\bibitem[AMo]{AbMo} S. Abhyankar, T.-T. Moh, {\em Embeddings of the line in the plane}, J. Reine Angew. Math. {\bf 276} (1975), 148-166.

\bibitem[AFKKZ]{AFKKZ} I.~V.~Arzhantsev, H.~Flenner, S.~Kaliman,
F.~Kutzschebauch, M.~Zaidenberg,
{\em Flexible varieties and automophism groups}.  Duke Math.\ J.\ {\bf 162} (2013), no. 4, 767--823.

 \bibitem[AZ]{AZ}  I. V. Arzhantsev, M. Zaidenberg, {\em Acyclic curves and group actions on affine toric surfaces.} Affine algebraic geometry, 1-41, World Sci. Publ., Hackensack, NJ, 2013.
 

\bibitem[AKuZ]{AKuZ} I. V. Arzhantsev, M. Zaidenberg, K. Kuyumzhiyan, {\em Flag varieties, toric varieties, and suspensions: three examples of infinite transitivity}, (Russian) 
Mat. Sb. {\bf 203} (2012), no. 7, 3-30; translation in Sb. Math. {\bf 203} (2012), no. 7-8, 923-949. 

\bibitem[AKuZ19]{AKuZ19} I. Arzhantsev, K.  Kuyumzhiyan, M. Zaidenberg, {\em Infinite transitivity, finite generation, and Demazure roots}, Adv. Math. {\bf 351} (2019), 1-32. 

\bibitem[BMS]{BMS} S. Bloch, M. Pavaman Murthy, L. Szpiro, {\em Zero cycles and the number of generators of an ideal}, {\bf 38}, 1989, Colloque en 
l\'honneur de Pierre Samuel (Orsay, 1987), pp. 51-74.





 




\bibitem[CLS]{CLS} Cox, David A.; Little, John B.; Schenck, Henry K. {\em Toric varieties}. Graduate Studies in Mathematics, {\bf 124}. 
American Mathematical Society, Providence, RI,  (2011) 841 pp. 

\bibitem[Cr]{Cr}  P. C. Craighero, {\em A result on m-flats in $\A_\bk^n$}, Rend. Sem. Mat. Univ. Padova {\bf 75} (1986), 39-46.


\bibitem[De]{De} M. Demazure,  {\em Sous-groupes alg\'ebriques de rang maximum du groupe de Cremona},
Annales scientifiques de l'\'Ecole Normale Sup\'erieure, S\'erie 4, Tome 3 (1970) no. 4,  507-588.



%


\bibitem[FS]{FS} P. Feller, I. van Santen, {\em Uniqueness of embeddings of the affine line into algebraic Groups}, 
 J. Algebraic Geom. {\bf 28} (2019), no. 4, 649-698.
 
 
\bibitem[FS21]{FS21} P. Feller, I. van Stampfli, {\em  Existence of embedding of smooth varieties into linear algebraic groups},  arXiv:2007.16164 (2020).

\bibitem[FKZ]{FKZ} H.~Flenner, S.~Kaliman, and M.~Zaidenberg, {\em A Gromov-Winkelmann type theorem for flexible varieties},  J. Eur. Math. Soc. (JEMS) 
{\bf18} (2016), no. 11, 2483-2510. 
%




\bibitem[Fr]{Fre} G. ~Freudenburg, {\em Algebraic Theory of Locally Nilpotent Derivations} , Encyclopaedia of Mathematical Sciences,  Springer,
 Berlin-Heidelberg-New York, 2006.




\bibitem[Gro58]{Gro58}  A.~ Grothendieck, {\em Torsion homologique et sections rationnelles}, Anneaux de Chow et Applications,
S\`eminaire Claude Chevalley, 1958, expos\^e n. 5, 1-29.



\bibitem[Ha]{Har} R.~Hartshorne, {\em Algebraic Geometry},
Springer-Verlag, New York-Heidelberg, 1977.

\bibitem[Hol]{Hol} Holme, Audun, {\em Embedding-obstruction for singular algebraic varieties in $\PP^N$}, Acta Math. {\bf 135} (1975), no. 3-4, 155-185.

 


\bibitem[Je]{Je}  Z. Jelonek, {\em The extension of regular and rational embeddings}, Math. Ann. {\bf 277} (1987), no. 1, 113-120.

\bibitem[Ka91]{Ka91} S.~ Kaliman, {\em Extensions of isomorphisms between affine algebraic subvarieties of $k^n$ to automorphisms of $k^n$}, 
Proc. Amer. Math. Soc. {\bf 113} (1991), no. 2, 325-334.




\bibitem[Ka20]{Ka20} S.~Kaliman, {\em Extensions of isomorphisms of subvarieties in flexible varieties}, Transform. Groups {\bf 25} (2020), no. 2, 517-575.

 



\bibitem[KaUd]{KaUd} S. ~Kaliman, D. ~Udumyan, {\em On automorphisms of flexible varieties}, arXiv:2008.02221 (2021).




\bibitem[Kl]{Kl} S.~L. Kleiman, {\em The transversality of a general translate},  Compositio Math. {\bf 28} (1974), 287-297. 

\bibitem[KR]{KrRu}  H.~Kraft, P.~ Russell, {\em Families of group actions, generic isotriviality, and linearization},
Transform. Groups 1{\bf 9} (2014), no. 3, 779-792.
 


\bibitem[Li]{Li} A. Liendo,  {\em Affine $\T$-varieties of complexity one and locally nilpotent derivations},
Transformation Groups volume {\bf 15} (2010), 389-425.


\bibitem[Lu]{Lu} D.~Luna, {\em Slices \'etales} (French) Sur les groupes alg\'ebriques, pp. 81-105. Bull. Soc. Math. France, Paris, M\'emoire {\bf 33} Soc. Math. France, Paris, 1973. 







\bibitem[PV]{PV}
V.~L.~Popov,  E.~B.~Vinberg, {\em Invariant Theory},  In Algebraic
geometry IV, ( A.\ N.\ Parshin, I.\ R.\ Shafarevich (eds.)),  Springer-Verlag, Berlin,
Heidelberg, New York, 1994.
%
%
%
\bibitem[Ra]{Ra} C.~P.~Ramanujam, {\em A note on automorphism groups of algebraic varieties}, Math.\ Ann.\
156 (1964), 25--33.







\bibitem[St]{St}  I. ~Stampfli, {\em Algebraic embeddings of $\C$ into $\SL_n(\C)$}, Transform. Groups {\bf 22} (2017), no. 2, 525-535.



\bibitem[Su]{Su} M. Suzuki, {\em Propi\'et\'es topologiques des polynomes de deux variables complexes, et automorphismes alg\'earigue de l'espace $\C^2$}, 
J. Math. Soc. Japan, {\bf  26} (1974), 241-257.

\bibitem[Sr]{Sr} V. Srinivas, {\em  On the embedding dimension of an affine variety},  Math. Ann., {\bf 289} (1991), no.1, 25-132.

\bibitem[Swan]{Swan} R. G. Swan, {\em A cancellation theorem for projective modules in the metastable range} Invent. Math. {\bf 27} (1974), 23-43.



\bibitem[Ud]{Ud} D. Udumyan, {\em Extension Problem for Flexible Varieties}, PhD thesis, Electronic Theses and Dissertations, University of Miami, 2019. 






\end{thebibliography}
\end{document}